\documentclass[11pt,reqno]{amsart}

\usepackage{ifdraft}
\usepackage{lmodern}
\usepackage[T1]{fontenc}

\usepackage{enumerate,tensor}
\usepackage{scalerel,mathtools}

\newtheorem{theorem}{Theorem}[section] % 1st argument is your name for it
\newtheorem{lemma}[theorem]{Lemma}     % 2nd argument is what is printed
\newtheorem{corollary}[theorem]{Corollary}

\theoremstyle{definition}
\newtheorem{defin}[theorem]{Definition}
\newtheorem{rem}[theorem]{Remark}
\newtheorem{exa}[theorem]{Example}

\numberwithin{equation}{section}

\usepackage{amsmath,latexsym,amssymb,amsxtra,cite}
\usepackage{amsfonts,epsfig,times,textcomp,upgreek,marginnote,mathtools,color}
\usepackage{enumerate,colortbl,setspace}

\usepackage[tight,scriptsize]{subfigure}
\usepackage{widetable,tikz,float}
\usetikzlibrary{arrows,automata,shapes,decorations}
\usetikzlibrary{decorations.pathmorphing,decorations.markings}
\usetikzlibrary{svg.path,calc}

\DeclareFontFamily{U}{mathx}{\hyphenchar\font45}
\DeclareFontShape{U}{mathx}{m}{n}{
      <5> <6> <7> <8> <9> <10>
      <10.95> <12> <14.4> <17.28> <20.74> <24.88>
      mathx10
      }{}
\DeclareSymbolFont{mathx}{U}{mathx}{m}{n}
\DeclareFontSubstitution{U}{mathx}{m}{n}
\DeclareMathAccent{\widebar}{0}{mathx}{"73}

\makeatletter
\def\hlinewd#1{
\noalign{\ifnum0=`}\fi\hrule \@height #1 
\futurelet\reserved@a\@xhline}
\makeatother

\newcommand{\ie}{\emph{i.e.}}
\newcommand{\resp}{\emph{resp.}}

\newcommand{\oo}{\hbox{$\mathbf\infty$}}
\newcommand{\DO}{\hbox{$\!\!\D^\frac{1}{3}\!\!$}}
\newcommand{\DT}{\hbox{$\!\!\D^\frac{2}{3}\!\!$}}
\newcommand{\DD}{\hbox{$\!\!\D^\frac{1}{1}\!\!$}}
\newcommand{\ccg}{\cellcolor{lightgray}}

\newcommand\dwideb[2]{%
  \setbox0=\hbox{$\widebar{#1}$}%
  \ht0=\dimexpr\ht0-.15ex\relax% CHANGE .15 TO AFFECT SPACING
  \widebar{\hspace*{-#2ex}\copy0\hspace*{#2ex}}%
}
\newcommand\dwb[2]{\hspace*{#2ex}\dwideb{#1}{#2}\hspace*{-#2ex}}%

\newcommand{\TT}[1]{\mathcal{T}_{#1}}
\newcommand{\TR}[1]{\mathcal{R}_{#1}}
\newcommand{\GT}[1]{G(\mathcal{T}_{#1})}
\newcommand{\MT}[1]{G^{+\!}(\mathcal{T}_{#1})}
\newcommand{\TTP}[1]{\mathcal{T}'_{#1}}
\newcommand{\GTP}[1]{G(\mathcal{T}'_{#1})}
\newcommand{\MTP}[1]{G^{+\!}(\mathcal{T}'_{#1})}
\newcommand{\VT}[1]{V(\mathcal{T}_{#1})}
\newcommand{\VTP}[1]{V(\mathcal{T}'_{#1})}

\tikzset{->-/.style={decoration={
  markings,
  mark=at position #1 with {\arrow{>}}},postaction={decorate}}}

\newcommand{\edg}[3]{\,\tikz[baseline=-.5ex]{\path[semithick] (0,0) edge	node[pos=0.25,above=-1,scale=.8] {$#2$}			%
													node[pos=0.75,above=-1,scale=.8] {$#3$} (#1ex,0);}\,}

\newcommand{\leaf}[1]{\,\tikz[baseline=-.5ex]{	\node [circle,draw=black,fill=lightgray,inner sep=1pt,minimum size=19pt,scale=.8] at (0,0) {$#1$};}\,}
\newcommand{\rootree}[8]{\,\tikz[baseline=-.5ex]{	\path[semithick] (0,0) edge[->-=.84,>=latex] 	node[pos=0.3,above,scale=.7] {$#3$}			%
																		node[pos=0.7,above,scale=.7] {$#4$} (#2ex,2ex);
										\path[semithick] (0,0) edge[->-=.84,>=latex]	node[pos=0.3,below,scale=.7] {$#5$}			%
																		node[pos=0.7,below,scale=.7] {$#6$} (#2ex,-2ex);
										\path (0,0) edge[white]	node[pos=0.7,scale=.7,black] {$\vdots$} (#2ex,0.6ex);
										\fill[left color=lightgray,right color=white] (#2ex,2ex) -- (#2ex+6ex,3.5ex) -- (#2ex+6ex,.5ex) -- cycle;
										\fill[left color=lightgray,right color=white] (#2ex,-2ex) -- (#2ex+6ex,-.5ex) -- (#2ex+6ex,-3.5ex) -- cycle;
										\node [circle,draw=black,fill=lightgray,inner sep=1pt,minimum size=19pt,scale=.85] at (0,0) {$#1$};
										\node [circle,draw=black,fill=lightgray,inner sep=1pt,minimum size=19pt,scale=.85] at (#2ex,2ex)	{$#7$};
										\node [circle,draw=black,fill=lightgray,inner sep=1pt,minimum size=19pt,scale=.85] at (#2ex,-2ex) {$#8$};}\,}
\newcommand{\staroot}[8]{\,\tikz[baseline=-.5ex]{	\path[semithick] (0,0) edge[->-=.87,>=latex]	node[pos=0.3,above,scale=.7] {$#3$}			%
																			node[pos=0.7,above,scale=.7] {$#4$} (#2ex,3ex);
										\path[semithick] (0,0) edge[->-=.87,>=latex]	node[pos=0.3,below,scale=.7] {$#5$}			%
																			node[pos=0.7,below,scale=.7] {$#6$} (#2ex,-3ex);
										\path (0,0) edge[white]	node[pos=0.7,scale=.7,black] {$\vdots$} (#2ex,0.6ex);
										\node [circle,draw=black,fill=lightgray,inner sep=1pt,minimum size=19pt,scale=.85] at (0,0) {$#1$};
										\node [circle,draw=black,fill=lightgray,inner sep=1pt,minimum size=19pt,scale=.85] at (#2ex,3ex)	{$#7$};
										\node [circle,draw=black,fill=lightgray,inner sep=1pt,minimum size=19pt,scale=.85] at (#2ex,-3ex) {$#8$};}\,}
\newcommand{\rootreeS}[1]{\,\tikz[baseline=-.5ex]{	\fill[left color=lightgray,right color=white] (0,0) -- (6ex,1.5ex) -- (6ex,-1.5ex) -- cycle;
										\node [circle,draw=black,fill=lightgray,inner sep=1pt,minimum size=19pt,scale=.85] at (0,0)	{$#1$};}\,}
\newcommand{\rootreeR}[1]{\,\tikz[baseline=-2ex]{	%\fill[left color=lightgray,right color=white] (0,0) -- (6ex,1.5ex) -- (6ex,-1.5ex) -- cycle;
										\fill[left color=lightgray,right color=white,scale=.85] (00ex+.25ex,-.25ex) -- (00ex+7ex,-3ex) -- (00ex+5ex,-6ex) -- cycle;
										\node [scale=.85] at (00ex+1.5ex,-4ex)	{$\TR{#1}$};
										\node [circle,draw=black,fill=lightgray,inner sep=1pt,minimum size=19pt,scale=.85] at (0,0)	{$#1$};}\,}
\newcommand{\rootreeRS}[1]{\,\tikz[baseline=-2ex]{	\path[semithick] (0,0) edge[->-=.84,>=latex]	node[pos=0.3,above,scale=.7] {$p_i$}			%->-=.84,>=latex
																			node[pos=0.7,above,scale=.7] {$q_i$} (09ex,0);
										\fill[left color=lightgray,right color=white,scale=.85] (00ex+.25ex,-.25ex) -- (00ex+7ex,-3ex) -- (00ex+5ex,-6ex) -- cycle;
										\fill[left color=lightgray,right color=white,scale=.85] (09ex,.25ex) -- (09ex+7ex,.5ex) -- (09ex+6ex,-4ex) -- cycle;
										\node [scale=.85] at (00ex+1.5ex,-4ex)	{$\TR{#1}$};
										\node [circle,draw=black,fill=lightgray,inner sep=1pt,minimum size=19pt,scale=.85] at (0,0)	{$#1$};
										\node [circle,draw=black,fill=lightgray,inner sep=1pt,minimum size=19pt,scale=.85] at (09ex,0)	{\phantom{$#1$}};}\,}
\newcommand{\rootreeH}[5]{\,\tikz[baseline=-.5ex]{	\path[semithick] (0,0) edge[->-=.84,>=latex]	node[pos=0.3,above,scale=.7] {$#3$}			%
																			node[pos=0.7,above,scale=.7] {$#4$} (#2ex,0);
										\node [circle,draw=black,fill=lightgray,inner sep=1pt,minimum size=19pt,scale=.85] at (0,0)	{$#1$};
										\node [circle,draw=black,fill=lightgray,inner sep=1pt,minimum size=19pt,scale=.85] at (#2ex,0)	{$#5$};}\,}
\newcommand{\vv}{v}

\newcommand{\vvx}{x}

\newcommand{\vvz}{z}

\newcommand{\MM}{M}
\newcommand{\GG}{G}
\newcommand{\BS}{\hbox{$\mathbf{B\!\!\;S}$}}

\newcommand{\ZZ}{{\mathbb Z}}
\newcommand{\BB}{{\mathbf B}}%braid
\newcommand{\Chiral}{\hbox{$M_\upchi$}}
\newcommand{\Knuth}{\hbox{$M_\upkappa$}}

\newcommand{\myroot}{h}
\newcommand{\myruut}{u}
\newcommand{\stable}{t}
\newcommand{\xx}{x}
\newcommand{\yy}{y}
\newcommand{\zz}{z}
\newcommand{\hh}{h}

\newcommand{\ee}{\upepsilon}
\newcommand{\mmu}{\upmu}
\newcommand{\autom}{\psi}
\newcommand{\inject}{\phi}
\newcommand{\dbar}{\mid\!\mid}
\newcommand{\length}[1]{\hbox{$\dbar\!\!#1\!\!\dbar$}}

\newcommand{\tta}{\mathtt{a}}
\newcommand{\ttb}{\mathtt{b}}
\newcommand{\ttc}{\mathtt{c}}

\newcommand{\tth}{\mathtt{h}}
\newcommand{\ttts}{\mathtt{s}}
\newcommand{\tttt}{\mathtt{t}}
\newcommand{\ttx}{\mathtt{x}}
\newcommand{\tty}{\mathtt{y}}
\newcommand{\ttz}{\mathtt{z}}

\newcommand{\D}{\Updelta}

\newcommand{\dR}{\mathbin{\backslash}}

\newcommand{\dL}{\mathbin{\slash}}
\newcommand{\jR}{\mathbin{\scriptstyle\vee}}
\newcommand{\jjR}{\mathbin{\bigvee}}
\newcommand{\QZ}{Q\!Z}
\newcommand{\gL}{\mathbin{\scriptstyle\wedge}}
\newcommand{\jL}{\mathbin{\widetilde{\scriptstyle\vee}}}
\newcommand{\gR}{\mathbin{\widetilde{\scriptstyle\wedge}}}

\newcommand\selecName{\uptheta}
\newcommand\selec[3]{\selecName_{#1}(#2,#3)}
\newcommand\onceMM{\widebar{M}}
\newcommand\onceGG{\widebar{G}}
\newcommand\seccMM{\dwb{M}{0.095}}
\newcommand\seccNu{\dwb{\upnu}{0.000}}
\newcommand\seccName{\dwb{\selecName}{0.031}}
\newcommand\secc[2]{\seccName(#1,#2)}%selec second

\newcommand\shuf[2]{{\overbracket[.5pt][1pt]{#1}}^{#2}}

\def\mx{\hbox{\color{black}{$\mmu_x$}}}
\def\my{\hbox{\color{black}{$\mmu_y$}}}
\def\mz{\hbox{\color{black}{$\mmu_z$}}}

\makeatletter
   \def\equalsfill{$\m@th\mathord=\mkern-7mu
     \cleaders\hbox{$\!\mathord=\!$}\hfill
     \mkern-7mu\mathord=$}
   \makeatother
\newcommand{\stack}[1]{{\ \stackrel{\text{\tiny #1}}{\hbox{\equalsfill}}\ }}

\usepackage{xcolor}
\definecolor{darkgreen}{rgb}{0.2,0.55,0}

\definecolor{orcidlogocol}{HTML}{A6CE39}
\tikzset{
  orcidlogo/.pic={
    \fill[orcidlogocol] svg{M256,128c0,70.7-57.3,128-128,128C57.3,256,0,198.7,0,128C0,57.3,57.3,0,128,0C198.7,0,256,57.3,256,128z};
    \fill[white] svg{M86.3,186.2H70.9V79.1h15.4v48.4V186.2z}
                 svg{M108.9,79.1h41.6c39.6,0,57,28.3,57,53.6c0,27.5-21.5,53.6-56.8,53.6h-41.8V79.1z M124.3,172.4h24.5c34.9,0,42.9-26.5,42.9-39.7c0-21.5-13.7-39.7-43.7-39.7h-23.7V172.4z}
                 svg{M88.7,56.8c0,5.5-4.5,10.1-10.1,10.1c-5.6,0-10.1-4.6-10.1-10.1c0-5.6,4.5-10.1,10.1-10.1C84.2,46.7,88.7,51.3,88.7,56.8z};
  }
}

\newcommand\orcidicon[1]{\href{https://orcid.org/#1}{\mbox{\scalerel*{
\begin{tikzpicture}[yscale=-1,transform shape]
\pic{orcidlogo};
\end{tikzpicture}
}{|}}}}

\usepackage{hyperref} %<--- Load after everything else

\begin{document}
\title[Cyclic amalgams, HNN extensions, and Garside one-relator groups]
{Cyclic amalgams, HNN extensions,\\ and Garside one-relator groups}
%{Tree products of cyclic groups, HNN extensions,\\ and Garside one-relator groups}
\author{Matthieu Picantin}% \orcidicon{0000-0002-7149-1770}
\address{IRIF, UMR 8243 CNRS \& Univ. Paris Diderot, 75013 Paris, France.}
\email{picantin@irif.fr}

%\remMP{orcid!!}
%\thanks{}

\subjclass[2020]{20E06 -- 20F05 -- 20F10 -- 20F36}

\keywords{amalgamated free product, HNN extension, tree product, knot group, braid monoid/group, one-relator group, Garside structure, automatic structure}

\date{}

\dedicatory{This paper is dedicated to the memory of my friend and mentor, Patrick Dehornoy.}

\begin{abstract}
Garside groups are a natural lattice-theoretic generalisation of the braid groups and spherical type Artin--Tits groups.
Here we show that the class of Garside groups is closed under some free products with cyclic amalgamated subgroups.
We deduce that every tree product of infinite cyclic groups is a Garside group.
Moreover, we study those cyclic HNN extensions of Garside groups that are Garside groups as well.
Using a theorem of~Pietrowski, we conclude this paper by stating that a non-cyclic one-relator group is Garside
if and only if its centre is nontrivial.
\end{abstract}

\maketitle

\section{Introduction}
\label{intro}

Braid groups are well understood due to Brieskorn's presentation theorem and the subsequent structural study by Deligne and Brieskorn--Saito~\cite{brieskornSaito,deligne}~:
their main combinatorial results express that every braid group is the group of fractions of a monoid in which divisibility has good properties,
and, in addition, there exists a distinguished element whose divisors encode the whole structure~: in modern terminology, such a monoid is called Garside.
The group of fractions of a Garside monoid is called a Garside group. Garside groups enjoy many remarkable group-theoretical, cohomological, and homotopy-theoretical properties \cite{dehornoyGarside,dehornoyTheory}.

\medskip
The aim of this paper is to explore further the class of Garside groups.
The latter happens to be closed under direct product and under some Zappa--Sz\'ep products.
We just recall that a \emph{Zappa--Sz\'ep product} (or \emph{bi-crossed product},
or \emph{knit product}) of groups or monoids~\cite{Szep,Zappa} is the natural extension of a semidirect product in which both groups (or monoids) act on one another,
the semidirect product corresponding to the case when one action is trivial, and the direct product to the case when both actions are trivial.
According to~\cite{dehornoyTheory}, the main general structural result about Garside structures known so far is that
every Garside group is an iterated Zappa--Sz\'ep product of Garside groups with infinite cyclic centre~\cite{pCenter}.

\medskip
Here we investigate the cyclic amalgamated free products of Garside monoids (Section~\ref{s:amalgam}).
From this we deduce that every tree product of infinite cyclic groups is a Garside group.
Then we characterise those cyclic HNN extensions of Garside groups that are again Garside groups (Section~\ref{s:hnnextension}).
As a nice consequence, using the solution of~Pietrowski for the isomorphism problem for one-relator groups with nontrivial centre~\cite{pietrowski},
we finally state that a non-cyclic one-relator group is Garside if and only if its centre is nontrivial (Section~\ref{s:pietrowski}).

%%%%%%%%%%%%%%%%%%%%%%%%%%
%%%%%%%%%%%%%%%%%%%%%%%%%%
%%%%%%%%%%%%%%%%%%%%%%%%%% 
\section{Background from Garside theory}\label{s:background}

In this section, we list some basic properties of
Garside monoids and groups.
For all the results and the examples quoted here, we refer the reader
to~\cite{dehornoyTheory}
(see also~\cite{dehornoyGarside,dehornoyComplete,dehornoySubword,dehornoyParis,pPhD,pCenter,pConjugacy,pExplicit,pTorus,pTransducer,pAutomaticon}).

\medbreak
Assume that~$M$ is a monoid. For~$a,b$ in~$M$,
we say that~$b$ is a left-divisor of~$a$---or that~$a$ is a right-multiple
of~$b$---if~$a=bd$ holds for some~$d$ in~$M$. An element~$c$ is a
lower common right-multiple---or a right-lcm---of~$a$ and~$b$
if it is a right-multiple of both~$a$ and~$b$, and every common right-multiple of~$a$
and~$b$ is a right-multiple of~$c$. Right-divisor, left-multiple, and
left-lcm are defined symmetrically.
For~$a,b$ in~$M$, we say that~$b$ \emph{divides}~$a$---or that~$b$ is a divisor of~$a$---if~$a=cbd$ holds for
some~$c,d$ in~$M$. 

\medbreak
The monoid~$M$ is said to be \emph{cancellative} when, for any~$a,b,c,d\in M$, $abc=adc$ implies~$b=d$.
And it is said to be \emph{conical} if~$1$ is its only invertible element, that is, $ab=1$ implies~$a=b=1$. 
Now, if~$c$, $c'$ are two right-lcms of~$a$ and~$b$,
necessarily~$c$ is a left-divisor of~$c'$, and~$c'$ is a left-divisor
of~$c$. If we assume~$M$ to be conical and cancellative, we have~$c=c'$.
In this case, the unique right-lcm of~$a$ and~$b$ is denoted by~$a
\jR b$, and the unique element~$c$ satisfying~$a
\jR b=ac$ is denoted by~$a\dR b$. We define the \emph{left-lcm}~$\jL$
and the left operation~$\dL$ symmetrically.  In particular, we
have\[a\jR b=a(a\dR b)=b(b\dR a),
\hbox{\quad  and\quad} a\jL b=(b\dL a)a=(a\dL b)b.\]Cancellativity and
conicity imply that left- and right-divisibility are order relations.

\begin{defin}\label{D:GarsideMonoid}
A monoid~$M$ is said to be \emph{Garside} if $M$ is conical and
cancellative, every pair of elements in~$M$ admits a left-lcm and a right-lcm,
and~$M$ admits a \emph{Garside element}, defined to be an element
whose left- and right-divisors coincide, are finite in number, and generate~$M$.
\end{defin}

\begin{exa}\label{e:classical} By~\cite{brieskornSaito}, all spherical type Artin--Tits monoids are Garside
monoids. The classical braid monoids of the complex reflection
groups~$G_{\tiny7},G_{\tiny11},G_{\tiny12},G_{\tiny13}, G_{\tiny15},G_{\tiny19},$ and~$G_{\tiny22}$ given
in~\cite{broueMalleRouquier} (see~\cite{dehornoyParis,pPhD}),
\hbox{Garside's hypercube monoids \cite{garside,pPhD}}, the dual braid
monoids for spherical type Artin--Tits groups \cite{bessisDual,bessisDigneMichel,birmanKoLee,pPhD,pExplicit}
and some so-called post-classical braid monoid~\cite{pBeer} for complex reflection groups of type~$(e,e,r)$ are also Garside monoids.
\end{exa}

\begin{defin}\label{d:envelop}The \emph{enveloping group} of a monoid~$\MM$
is defined as the group admitting the same presentation as~$\MM$;
formally, it can be viewed as the monoid~$\langle~{\MM}\cup{\overline\MM}:{R}_\MM\cup{F}_\MM~\rangle^+$,
where~$R_\MM$ is the family of all relations~$bc = d$ with~$b,c,d\in\MM$
and~$F_\MM$ is the family of all relations~$b\overline{b}=\overline{b}b=1$ with~$b\in\MM$.
\end{defin}

If~$M$ is a Garside monoid, then~$M$ satisfies Ore's
conditions~\cite{cliffordPreston}, and it embeds in a group of right-fractions, and,
symmetrically, in a group of left-fractions. In this case, by
construction, every right-fraction~$ab^{-1}$ with~$a,b$ in~$M$ can be
expressed as a left-fraction~$c^{-1}d$ with~$c,d$ in~$M$, and conversely. Therefore, the two
groups coincide, and  there is no ambiguity in speaking of \emph{the}
group of fractions of a Garside monoid.

\begin{defin}\label{D:garsideGroup}
A group~$G$ is said to be a \emph{Garside group}
if there exists a Garside monoid of which~$G$ is the enveloping group, hence the group of fractions.
\end{defin}

\begin{exa}\label{e:Chiral}
Let us consider the monoid~$\Chiral$ with presentation
\[\langle~\tta,\ttb,\ttc: \tta^2 = \ttb\ttc, \ttb^3 = \ttc\tta~\rangle^{+}.\]
The monoid~$\Chiral$ is a typical example of a Garside monoid, and, in
addition, $\Chiral$ has the distinguishing feature to not be
anti-automorphic, contrary to those examples mentioned in~Example~\ref{e:classical}.
Its group of fractions happens to be the group~$\langle~\tta,\ttb:\tta^3=\ttb^4~\rangle$ of the $(3,4)$-torus knot. We refer to~\cite{pPhD,pTorus} for more Garside structures for torus link groups.
\end{exa}

\begin{rem}Garside structures had been called \emph{small Gaussian} or \emph{thin Gaussian} in earlier works until~2000 \cite{dehornoyParis,pPhD,pCenter,pConjugacy}.
\end{rem}

\begin{lemma}\label{lem-calculous}{\rm\cite{dehornoyTheory}}
Assume that~$M$ is a Garside monoid. Then the following identities hold
in~$M$:
\begin{gather*}
(ab)\jR(ac)=a(b\jR c),\\
c\dR(ab)=(c\dR a)((a\dR c)\dR b),\hbox{\qquad} (ab)\dR c =b\dR(a\dR c),\\
(a\jR b)\dR c =(a\dR b)\dR(a\dR c) =(b\dR a)\dR(b\dR c),
\hbox{\qquad} c\dR(a\jR b)=(c\dR a)\jR(c\dR b).
\end{gather*}
\end{lemma}

\begin{lemma}\label{lem-atomic}{\rm\cite{dehornoyTheory}}
Assume that~$M$ is a Garside monoid. The following equivalent
assertions hold~:

 (i) There exists a mapping~$\upnu$ from~$M$ to the integers
satisfying~$\upnu(a)>0$ for every~${a\not=1}$ in~$M$, and
%satisfying
~$\upnu(ab)\geq\upnu(a)+\upnu(b)$ for any~$a,b$ in~$M$;

 (ii) For every set~$X$ that generates~$M$ and for every~$a$ in~$M$, the
lengths of the decompositions of~$a$ as products of elements in~$X$ have
a finite upper bound.
\end{lemma}

\begin{defin}\label{D:Norm}
A monoid is said to be \emph{atomic} if it satisfies the equivalent
conditions of Lemma~\ref{lem-atomic}. An \emph{atom} is defined
to be a nontrivial element~$a$ such that~$a=bc$ implies~$b=1$ or~$c=1$.
The \emph{norm} function~\hbox{$\length{.}$} of an atomic monoid~$M$ is defined
in such a way that, for every~$a$ in~$M$,
$\length{a}$ is the upper bound of the lengths of the decompositions of~$a$
as products of atoms.
\end{defin}

\begin{exa}\label{ex-knuth}
The monoid~$\Knuth$ defined by the presentation\[\langle~\ttx,\tty:
\ttx\tty\ttx\tty\ttx=\tty\tty~\rangle^{+}\] is another example of a Garside monoid, which 
%as for it,
admits no additive norm, \ie, no norm~$\upnu$
satisfying~$\upnu(ab)=\upnu(a)+\upnu(b)$ for any two elements~$a,b$ in~$\Knuth$.
Its group of fractions happens to be isomorphic to~$\langle~\tty,\ttz:\tty^3=\ttz^3~\rangle$.
See~\cite{pPhD} for further details about~$\Knuth$.
\end{exa}

By the previous lemma, every element in a Garside
monoid has finitely many left-divisors, only then, for every pair of
elements~$(a, b)$, the common left-divisors of~$a$ and~$b$ admit a right-lcm,
which is therefore the \emph{left-gcd} of~$a$ and~$b$.
This left-gcd will be denoted by~$a\gL b$.
We define the \emph{right-gcd}~$\gR$ symmetrically. 

\begin{lemma}\label{T:Closure}{\rm\cite{dehornoyTheory}}
Assume that~$M$ is a Garside monoid. Then  it admits a finite
generating subset that is closed under~$\dR,\dL,\jR,\gL,\jL,$ and~$\gR$.
\end{lemma}

Every Garside monoid admits a finite set of atoms, and
this set is the minimal generating set~\cite{dehornoyParis}. The hypothesis that
there exists a finite generating subset that is closed under~$\dR$
implies that the closure of the atoms under~$\dR$ is finite.

\begin{defin}\label{def-GarsideElement} If~$M$ is a Garside monoid, the
closure~$S$ of its atoms under~$\dR$ and~$\jR$ is finite---its elements are
called \emph{simple elements}, and their right-lcm is the (minimal) Garside element denoted by~$\D$.
The finite lattice~$(S,\gL,\jR,1,\D)$ nicely summarises the whole structure of~$M$ in a nutshell.
\end{defin}

\begin{defin}\label{def-exponent} The set of all Garside elements of~$M$ is~$\D^+=\{\D^p:p>0\}$ (see~\cite{pCenter} for instance)
and we denote by~$\D^{\frac{p}{q}}$ the set of the $q$-th roots of the element~$\D^p$.
For any root~$\myroot$ of a Garside element,
we denote by~$(\updelta(\myroot),\ee(\myroot))$ the lexicographically minimal pair~$(p,q)$
satisfying~$\myroot\in\D^{\frac{ep}{q}}$, 
where~$e$ is the smallest positive integer for which~$\D^e$ is central, hence~$e=\ee(\D)$.
In particular, the so-called \emph{exponent}~$\ee(\myroot)$ of such a root~$\myroot$
is the smallest integer~$q$ such that~$\myroot^q$ is central.
\end{defin}

\begin{exa}\label{E:Lattice}
The lattices of simple elements in~$\Chiral$ of~Example~\ref{e:Chiral}
and in~$\Knuth$ of~Example~\ref{ex-knuth} are displayed
in~Figures~\ref{f:chi9} and~\ref{f:knuth}.
More lattices are displayed in Figures~\ref{fig-amalgam}, \ref{fig-hnnKnuth96}, and~\ref{fig-pietrow}.
Examples~\ref{ex-amalOK} and~\ref{ex-hnnKnuth} will include the computation of some small sets~$\D^{\frac{p}{q}}$ for~$\Chiral$ and for~$\Knuth$.
\end{exa}

\begin{figure}[!t]
\centering
\subfigure[\label{f:chi9}The 9-simple lattice of~$\Chiral$.]{
	\hspace*{17mm}
\scalebox{.65}{
\begin{tikzpicture}[yscale=.9,scale=.9,line width=1pt]
  \node (0) at (71.0bp,7.0bp) [draw,ellipse] {};
  \node (a) at (39.0bp,57.0bp) [draw,ellipse] {};
  \node (c) at (106.0bp,107.0bp) [draw,ellipse] {};
  \node (b) at (71.0bp,57.0bp) [draw,ellipse] {};
  \node (aaa) at (42.0bp,207.0bp) [draw,ellipse] {};
  \node (bb) at (74.0bp,107.0bp) [draw,ellipse] {};
  \node (ca) at (77.0bp,157.0bp) [draw,ellipse] {};
  \node (aa) at (42.0bp,107.0bp) [draw,ellipse] {};
  \node (ab) at (7.0bp,157.0bp) [draw,ellipse] {};
\path [orange]			(0) edge node[left=-1,pos=.47]{$\tta$}		(a);
\path [blue,dashed]		(0) edge node[left=-2]{$\ttb$}			(b);
\path [darkgray,dotted]	(0) edge node[right=-1,pos=.205]{$\ttc$}	(c);
\path [orange] (a) edge (aa);
\path [blue,dashed] (ca) edge (aaa);
\path [darkgray,dotted] (b) edge (aa);
\path [blue,dashed] (bb) edge (ca);
\path [blue,dashed] (b) edge (bb);
\path [orange] (aa) edge (aaa);
\path [blue,dashed] (a) edge (ab);
\path [orange] (c) edge (ca);
\path [darkgray,dotted] (ab) edge (aaa);
  \node (aaa) at (45.0bp,222.0bp) [] {$\D_\upchi$};
\end{tikzpicture}
}
	\hspace*{17mm}
}
\subfigure[\label{f:knuth}The 12-simple lattice of~$\Knuth$.]{
	\hspace*{17mm}
\scalebox{.65}{
\begin{tikzpicture}[yscale=.9,scale=.9,line width=1pt]
  \node (0) at (25.0bp,7.0bp) [draw,ellipse] {};
  \node (a) at (9.0bp,57.0bp) [draw,ellipse] {};
  \node (b) at (42.0bp,57.0bp) [draw,ellipse] {};
  \node (ba) at (77.0bp,107.0bp) [draw,ellipse] {};
  \node (bb) at (42.0bp,257.0bp) [draw,ellipse] {};
  \node (aba) at (7.0bp,157.0bp) [draw,ellipse] {};
  \node (babab) at (75.0bp,257.0bp) [draw,ellipse] {};
  \node (bbb) at (58.0bp,307.0bp) [draw,ellipse] {};
  \node (baba) at (77.0bp,207.0bp) [draw,ellipse] {};
  \node (ab) at (7.0bp,107.0bp) [draw,ellipse] {};
  \node (bab) at (77.0bp,157.0bp) [draw,ellipse] {};
  \node (abab) at (7.0bp,207.0bp) [draw,ellipse] {};
\path [orange, densely dashed]	(0) edge node[left]{$\ttx$}		(a);
\path [blue, dotted]	(0) edge node[right,pos=.465]{$\tty$}	(b);
\draw [blue, dotted] (bb) ..controls (47.033bp,273.1bp) and (52.941bp,290.82bp)  .. (bbb);
\draw [orange, densely dashed] (ab) ..controls (7.0bp,123.53bp) and (7.0bp,140.71bp)  .. (aba);
\draw [blue, dotted] (aba) ..controls (7.0bp,173.53bp) and (7.0bp,190.71bp)  .. (abab);
\draw [orange, densely dashed] (bab) ..controls (77.0bp,173.53bp) and (77.0bp,190.71bp)  .. (baba);
\draw [blue, dotted] (baba) ..controls (76.353bp,223.53bp) and (75.637bp,240.71bp)  .. (babab);
\draw [orange, densely dashed] (abab) ..controls (17.447bp,222.33bp) and (31.417bp,241.49bp)  .. (bb);
\draw [blue, dotted] (ba) ..controls (77.0bp,123.53bp) and (77.0bp,140.71bp)  .. (bab);
\draw [blue, dotted] (b) ..controls (42.0bp,93.647bp) and (42.0bp,220.58bp)  .. (bb);
\draw [orange, densely dashed] (babab) ..controls (69.653bp,273.1bp) and (63.375bp,290.82bp)  .. (bbb);
\draw [blue, dotted] (a) ..controls (8.3531bp,73.525bp) and (7.637bp,90.712bp)  .. (ab);
\draw [orange, densely dashed] (b) ..controls (52.447bp,72.327bp) and (66.417bp,91.486bp)  .. (ba);
  \node (aaa) at (60.0bp,322.0bp) [] {$\D_\upkappa$};
\end{tikzpicture}
}
	\hspace*{17mm}
}\caption{Two examples of lattices of simples from~Examples~\ref{e:Chiral} and~\ref{ex-knuth}.}
\end{figure}

\medskip
Any Garside group can be the enveloping group of various monoids,
many of whom can be Garside monoids.
Now its centre gives rise to common structural constraints on all of them. 

\begin{lemma}\label{lem-center-quotient}{\rm{\cite{pCenter}}} The centre of the group of fractions of a Garside monoid~$\MM$ is the group of fractions of the centre of~$\MM$.
\end{lemma}

\begin{theorem}\label{thm-quasicenter}{\rm{\cite{pCenter}}} The quasi-centre~$\QZ$ (\resp\ the centre) of a Garside monoid~$\MM$ is a free abelian submonoid of~$\MM$,
and the function~$a\mapsto\D_a=\jjR\,M\!\dR\!a$ is a surjective semilattice homomorphism from~$(\MM,\jR)$ onto~$(\QZ,\jR)$.
\end{theorem}

\begin{corollary}\label{cor-quasicenter} The free abelian group of rank~$n$ is the group of fractions
of a unique Garside monoid up to isomorphism.
\end{corollary}

\medskip
We conclude this section by recalling how to
effectively recognise Garside monoids.

\smallbreak
\begin{defin}\label{def-complemented}
A monoid presentation~$\langle~A:R~\rangle^{+}$ is called \emph{right-complemented}
if $R$ contains no $\varepsilon$-relation (that is, no relation~$w = \varepsilon$ with~$w$ nonempty),
no relation $s\cdots = s\cdots$ with~$s\in{A}$ and, for~$s\not=t\in{A}$, at most one relation~$s\cdots = t\cdots$.
\end{defin}

\smallbreak
\begin{defin}\label{def-syntactic}
A \emph{syntactic right-complement} on an alphabet~$A$ is a partial map~$\selecName$ from~$A^2$ to~$A^*$
such that~$\selecName(x,x) = \varepsilon$ holds for every~$x\in{A}$ and, if~$\selecName(x,y)$ is defined, then so is~$\selecName(y,x)$.
\end{defin}

\begin{lemma}\label{lem-complem-selec}{\rm{\cite{dehornoyTheory}}} %[Lemma 4.3]
A monoid presentation~$\langle~A:R~\rangle^{+}$ is right-complemented
if and only if there exists a syntactic right-complement~$\selecName$
such that~$R$ consists of all relations~$x\selecName(x,y) = y\selecName(y,x)$ with~$x\not=y\in{A}$.
\end{lemma}

In the situation of Lemma~\ref{lem-complem-selec}, we naturally say that the presentation~$\langle~A:R~\rangle^{+}$
is associated with the syntactic right-complement~$\selecName$ (which is uniquely determined by the presentation),
and then we write~$\langle~A:R_{\selecName}~\rangle^{+}$.

\begin{lemma}\label{lem-selec}{\rm{\cite{dehornoyTheory}}} %[Lemma 4.6]
Assume that~$\langle~A:R_{\selecName}~\rangle^{+}$ is a right-complemented presentation.
Then there exists a unique minimal extension of the 
syntactic right-complement~$\selecName$
into a partial map~$\selecName$ from~$A^*\times A^*$ to~$A^*$ that satisfies the rules
\begin{equation} \selec{}{u}{\varepsilon}=\varepsilon,	\quad\selec{}{\varepsilon}{u}=u,	\quad\selec{}{uv}{uw}=\selec{}{v}{w},		\label{selec-epsilon}		\tag{$\selecName$-epsilon}		\end{equation}
\begin{equation} \selec{}{u}{vw}=\selec{}{u}{v}\ \selec{}{\selec{}{v}{u}}{w},
	\hbox{\quad and\quad}	\selec{}{vw}{u}=\selec{}{w}{\selec{}{v}{u}}.											\label{selec-extension}	\tag{$\selecName$-extension}		\end{equation}
The map~$\selecName$ is such that~$\selec{}{u}{v}$ exists if and only if $\selec{}{v}{u}$ does.
\end{lemma}

\begin{defin}\label{def-cube}
Assume that $\langle~A:R_{\selecName}~\rangle^{+}$ is a right-complemented presentation. 
A triple~$(u,v,w)$ of words in~$A^*$ satisfies the \emph{$\selecName$-cube condition} whenever it satisfies
\begin{equation} \selec{}{\selec{}{u}{v}}{\selec{}{u}{w}}=\selec{}{\selec{}{v}{u}}{\selec{}{v}{w}},\label{selec-cube}	\tag{$\selecName$-cube}		\end{equation}
meaning that either both sides are defined and they are equal, or neither is defined (see~Figure~\ref{fig-cube}).
The $\selecName$-cube condition is satisfied on~$S\subseteq A^*$ if every triple of words in~$S$ satisfies~it.
\begin{figure}[!ht]
\begin{center}
\begin{tikzpicture}[thick,scale=2.49,minimum size=0pt,inner sep=1pt]
	\node (00) at (0,0) {};
	\node (11) at (1,1) {};
	\node (31) at (3,1) {};
	\node (02) at (0,2) {};
	\node (22) at (2,2) {};
	\node (13) at (1,3) {};
	\node (33) at (3,3) {};
	\node (20a) at (1.7,0.15) {};
	\node (20b) at (2.3,0.3) {};
	\node (20c) at (2.1,-0.35) {};
	\node[inner sep=1pt] (20aa) at (20a) {};
	\node[inner sep=1pt] (20bb) at (20b) {};
	\node[inner sep=1pt] (20cc) at (20c) {};
	\path[->,>=latex]
		(11)	edge 		node[above=2,near start=5]{$\selec{}{v}{w}$}			(31);
	\fill[color=lightgray] (0,0) -- (1.7,0.15) -- (2.1,-0.35) -- (0,0);	
	\fill[color=lightgray] (2,2) -- (1.7,0.15) -- (2.3,0.3) -- (2,2);	
	\fill[color=lightgray] (3,1) -- (2.3,0.3) -- (2.1,-0.35) -- (3,1);	
	\draw[color=white,line width=5pt] (0,0) -- (1.7,0.15);
	\draw[color=white,line width=5pt] (0,0) -- (2.1,-0.35);	
	\draw[color=white,line width=5pt] (1.7,0.15) -- (2.1,-0.35);	
	\draw[color=white,line width=5pt] (3,1) -- (2.3,0.3);	
	\draw[color=white,line width=5pt] (3,1) -- (2.1,-0.35);	
	\draw[color=white,line width=5pt] (2.3,0.3) -- (2.1,-0.35);	
	\draw[color=white,line width=5pt] (2,2) -- (1.7,0.15);	
	\draw[color=white,line width=5pt] (2,2) -- (2.3,0.3);	
	\draw[color=white,line width=5pt] (1.7,0.15) -- (2.3,0.3);	
	\path[->,>=latex]
		(13)	edge 		node[left=4,near start]{$u$}		(02)
		(13)	edge 		node[above=2]{$w$}				(33)
		(13)	edge 		node[right=2,near start]{$v$}		(11)
		(02)	edge 		node[left]{$\selec{}{u}{v}$}		(00)
		(02)	edge 		node[below=2,near start]{$\selec{}{u}{w}$}	(22)
		(33)	edge 		node[above,sloped]{$\selec{}{w}{u}$}						(22)
		(33)	edge 		node[right]{$\selec{}{w}{v}$}			(31)
		(11)	edge 		node[right=4,pos=.45]{$\selec{}{v}{u}$}					(00)
		(00)	edge[gray] 	node[above,sloped,pos=.55]{$\selec{}{\selec{}{u}{v}}{\selec{}{u}{w}}$}	(20a)
		(00)	edge[gray] 	node[below,sloped,pos=.43]{$\selec{}{\selec{}{v}{u}}{\selec{}{v}{w}}$}	(20c)
		(22)	edge[gray] 	node[above,near start=5]{}			(20a)
		(22)	edge[gray] 	node[above,near start=5]{}			(20b)
		(31)	edge[gray] 	node[above,near start=5]{}			(20b)
		(31)	edge[gray] 	node[above,near start=5]{}			(20c)
		;
\end{tikzpicture}
\caption{The $\selecName$-cube condition: when one can draw the six faces of the cube,
then each of the three small gray triangular sectors is labelled by two equal words, and the cube closes. See Definition~\ref{def-cube}.}
\label{fig-cube}
\end{center}
\end{figure}
\end{defin}

\medbreak Rules~\eqref{selec-epsilon} and~\eqref{selec-extension}, and Condition~\eqref{selec-cube} on words have to be compared
with those of~Lemma~\ref{lem-calculous} on elements:
any Garside monoid with set of atoms~$A$ admits a right-complemented presentation~$\langle~A:R_{\selecName}~\rangle^{+}$
satisfying Condition~\eqref{selec-cube} on~$A^*$.

\begin{exa}\label{ex-complem}
The monoid~$\Chiral$ from~Example~\ref{e:Chiral} admits the right-complemented presentation~$\langle~A:R_{\selecName_\upchi}~\rangle^{+}$
with~$\selec{\upchi}{\tta}{\ttb}=\tta$, $\selec{\upchi}{\tta}{\ttc}=\tta^2$,
$\selec{\upchi}{\ttb}{\tta}=\ttc$, $\selec{\upchi}{\ttb}{\ttc}=\ttb^2$,
$\selec{\upchi}{\ttc}{\tta}=\tta\ttb$, and~$\selec{\upchi}{\ttc}{\ttb}=\tta$
(one could choose~$\selecName'_\upchi$ with~$\selecName'_\upchi(\tta,\ttc)=\ttb\ttc$ instead).
One can compute for instance~$\selec{\upchi}{\selec{\upchi}{\tta}{\ttb}}{\selec{\upchi}{\tta}{\ttc\ttb}}=\selec{\upchi}{\tta}{\tta\tta\ttc}=\tta\ttc=\selec{\upchi}{\ttc}{\ttb\ttb\tta}
=\selec{\upchi}{\selec{\upchi}{\ttb}{\tta}}{\selec{\upchi}{\ttb}{\ttc\ttb}}$, witnessing the $\selecName_\upchi$-cube condition for the triple~$(\tta,\ttb,\ttc\ttb)$.
\end{exa}

\medbreak The criterion we shall use in the sequel is~:

\begin{theorem}\label{thm-criterion}{\rm{\cite{dehornoyTheory}}} %[Proposition 4.16]
Assume that an atomic monoid~$M$ admits
a right-complemented presentation~$\langle~A:R_{\selecName}~\rangle^{+}$
satisfying the $\selecName$-cube condition on~$A$.
Then $M$ is left-cancellative and admits conditional right-lcms, that is, any two elements of~$M$
that admit a common right-multiple admit a right-lcm.
\end{theorem}

Moreover, if~$M$ is also right-cancellative and admits a Garside element, then $M$ is a Garside monoid.
For alternative Garsidity criteria and details, we refer to~\cite{dehornoyTheory}
(see also~\cite{dehornoyGarside,dehornoyComplete,dehornoySubword,dehornoyParis}).

%%%%%%%%%%%%%%%%%%%%%%%%%%
%%%%%%%%%%%%%%%%%%%%%%%%%%
%%%%%%%%%%%%%%%%%%%%%%%%%% 
\section{Amalgamated free products}\label{s:amalgam}

We prove that the class of Garside groups is closed under some free products with cyclic amalgamated subgroups.
This turns out to be exactly what we need to deduce that every tree product of infinite cyclic groups is a Garside group.

\medskip
\begin{defin}\label{d:amalgam} Let~$\MM_1$, $\MM_2$, $H$ be 
monoids with morphisms~$\inject_1:H\hookrightarrow\MM_1$
and~$\inject_2:H\hookrightarrow\MM_2$. 
The \emph{amalgamated free product} of~$\MM_1$ and~$\MM_2$ with respect to~$H$, $\inject_1$, and~$\inject_2$ is the monoid
\[\langle~\MM_1\star\MM_2~: \inject_1(h)=\inject_2(h), h\in H~\rangle^{+}.\]
When $H =\langle~\myroot~\rangle^+$ is cyclic, we simply write~$\inject_1(\myroot) =\myroot_1$, $\inject_2(\myroot) =\myroot_2$,
and the amalgamated free product is denoted by~$\MM_1\star_{\myroot_1=\myroot_2}\MM_2$.
\end{defin}
 
\begin{theorem}\label{thm-amalGarside}
Let~$\MM_1$ and~$\MM_2$ be some Garside monoids. Then, for any root~$\myroot_1$ of any Garside element in~$\MM_1$
and  any root~$\myroot_2$ of any Garside element in~$\MM_2$,
the cyclic amalgamated free product~$\MM_1\star_{\myroot_1=\myroot_2}\MM_2$ is a Garside~monoid.
\end{theorem}

Actually, a necessary assumption is that $\inject_i(H)$ has to contain a Garside element of~$\MM_i$ for~$i\in\{1,2\}$.
When restricted to cyclic amalgamated submonoids, the latter naturally expresses in terms of roots of Garside elements.

\begin{figure}[!ht]
\centering
\vspace*{-20pt}
	\scalebox{.96}{
	\begin{tikzpicture}[>=latex,very thick,top color=white,bottom color=gray]
	\tikzstyle{every state}=[minimum size=6pt,inner sep=1pt]
	\node (m3) at (9,8.5) {~};
	\clip[decorate, decoration={random steps,segment length=2mm}] (-1,-.6) rectangle (11,7.5);
	\node (1) at (5,0) {};	
	\node (a1) at (5,2) {};
	\node (a2) at (5,4) {};
	\node (a3) at (5,6) {};
	\node (a4) at (5,8) {};
	\node (m1) at (1,7.5) {};
	\node (m2) at (9,7.5) {};
	\def\pathone{(m1) 	to[out=270, in=180] (1.center) to[bend left]  (a1.center)
					to[bend left] (a2.center) to[bend left] (a3.center) to[bend left] (a4.center)};
	\def\pathtwo{(m2) 	to[out=270, in=0] (1.center) to[bend right]  (a1.center)
					to[bend right] (a2.center) to[bend right] (a3.center) to[bend right] (a4.center)};
	\fill[bottom color=orange!30, top color=white] \pathone;
	\fill[bottom color=blue!30, top color=white] \pathtwo;
	\path (1) 	edge[orange, bend left]  (a1.center)
		(a1.center) edge[orange, bend left] (a2.center)
		(a2.center) edge[orange, bend left] (a3.center)
		(a3.center) edge[orange, bend left] (a4.center);
	\path (1) 	edge[blue, bend right] (a1.center)
		(a1.center) edge[blue, bend right] (a2.center)
		(a2.center) edge[blue, bend right] (a3.center)
		(a3.center) edge[blue, bend right] (a4.center);
	\node[state, fill=white, label=west:$b_1$] (b1) at (3.4,2.3) {};
	\node[state, fill=white, label=east:$b_2$] (b2) at (7,3) {};
	\node[state, fill=white,minimum size=4pt,fill=white] (b1s) at (4.73,3.3) {};
	\node[state, fill=white,minimum size=4pt,fill=white] (b2s) at (5.29,5) {};
	\path	(1) edge[orange, bend left,  decorate, decoration={random steps,segment length=2mm}] (b1)
		(b1) edge[orange, bend left,  decorate, decoration={random steps,segment length=2mm}] (b1s)
		(1) edge[blue, bend right,  decorate, decoration={random steps,segment length=2mm}] (b2)
		(b2) edge[blue,  decorate, decoration={random steps,segment length=2mm}] (b2s);
	\node[state, fill=white, label=south east:$1$] (1) at (5,0) {};	
	\node[state, fill=white, label=west:$\myroot_1$, label=east:$\myroot_2$] (a1) at (5,2) {};
	\node[state, fill=white, label=west:$\myroot_1^2$, label=east:$\myroot_2^2$] (a2) at (5,4) {};
	\node[state, fill=white, label=west:$\myroot_1^3$, label=east:$\myroot_2^3$] (a3) at (5,6) {};
	\node (mb1) at (3,3.95) {$\mmu_{b_1}=2^{~}$};
	\node (mb2) at (7,5.98) {$\mmu_{b_2}=3^{~}$};
	\node[label=south east:$M_1$] (m1) at (1,7.5) {};
	\node[label=south west:$M_2$] (m2) at (9,7.5) {};
\end{tikzpicture}}
\vspace*{-10pt}
\caption{The lattice structure of the restriction to~$M_1\jR M_2$ of the amalgamated free product~$M_1\star_{\myroot_1=\myroot_2}M_2$.
See~the proof of Theorem~\ref{thm-amalGarside}.}
\label{f:amalgam}
\end{figure}

\begin{rem}\label{rem-extraction} Some algorithms for root extraction in~Garside groups have been proposed in~\cite{leeExtraction,sibert,styshnev}~:
the extraction problem of an $n$-th root in a Garside group~$G$ reduces to a conjugacy problem in the Garside group~$\ZZ\ltimes G^n$
and hence is decidable (see~\cite{pCenter,pConjugacy}).
Let us mention that the number of integers~$n$ for which an element admits an $n$-th root is finite
and that the number of conjugacy classes of the $n$-th roots of an element is finite.
\end{rem}

\begin{proof}[Proof of Theorem~\ref{thm-amalGarside}]
First, by hypothesis, the monoids~$\MM_1$ and~$\MM_2$ are cancellative 
and the amalgamated monoid~$H$ is the (infinite) cyclic monoid,
hence $\inject_i(H)=\langle~\myroot_i~\rangle^{+}$ is a so-called \emph{unitary} submonoid of~$\MM_i$ for each~$i\in\{1,2\}$,
that is, either~$hb\in\inject_i(H)$ or~$bh\in\inject_i(H)$ together with~$h\in\inject_i(H)$ imply~$b\in\inject_i(H)$.
According to~\cite[Corollary~3.4]{howie62} (see also~\cite{dekov91}),
the monoids~$\MM_1$ and~$\MM_2$ both embed into the cyclic amalgamated product~$\MM_1\star_{\myroot_1=\myroot_2}\MM_2$.

\medskip
Next, the root assumption guarantees~$\MM$ to inherit the atomicity from~$\MM_1$ and~$\MM_2$.
Indeed, the lengths of the decompositions of any element~$b\in\MM$ as products of atoms 
can be shown to be 
upper-bounded by the length of some central Garside element common to~$\MM_1$, to~$\MM_2$, and to~$H$.
Formally, by denoting~$S_j=\MM_j\setminus \myroot_j\MM_j$ for~$j\in\{1,2\}$,
it is known that the natural map~$\psi_j: \inject_j(H)\times S_j\rightarrow\MM_j$ is bijective (see~\cite{dekov98} for instance).
We deduce that any element~$b\in\MM$ admits a unique decomposition~$\myroot_1^pb_1b_2\cdots b_{2n}$ with~$b_{2k-1}\in S_1$, $b_{2k}\in S_2$, and~$1\leq k\leq n$.
Now, for~$j\in\{1,2\}$, each element~$b_{2k-j}\in S_j$ admits a minimal
right-multiple of the form~$b_{2k-j}b'_{2k-j}=\myroot_j^{\ee(\myroot_j)f_{2k-j}}$ for some~$b'_{2k-j}\in\MM_j$ and some integer~$f_{2k-j}>0$.
Therefore we obtain
\[b\,{b}'_{2n}\cdots{b}'_{2}{b}'_{1}=\myroot_1^p\myroot_2^{\ee(\myroot_2)f_{2n}}\cdots\myroot_2^{\ee(\myroot_2)f_{2}}\myroot_1^{\ee(\myroot_1)f_{1}}
=\myroot_1^{N_{b}}=\myroot_2^{N_{b}},\]
where~$N_{b}$ is the uniquely determined number~$p+\ee(\myroot_1)\sum_{k=1}^{n}f_{2k-1}+\ee(\myroot_2)\sum_{k=1}^{n}f_{2k}$.
Denoting by~$\length{\cdot}_1,\length{\cdot}_2$, and~$\length{\cdot}$ the norms of~$\MM_1,\MM_2$, and~$\MM$, respectively, we conclude
\[\length{b}\ \leq\max(
\length{\myroot_1^{N_{b}}}_1,
\length{\myroot_2^{N_{b}}}_2),\]
which gives the claim.

\medskip
Let~$A_i$ be the set of atoms of~$\MM_i$ and let~$\selecName_i$ be a syntactic right-complement
such that~$\MM_i$ admits the complemented presentation~$\langle~A_i:R_{\selecName_i}~\rangle^{+}$ for~$i\in\{1,2\}$.
By definition, the cyclic amalgamated free product monoid~$\MM=\MM_1\star_{\myroot_1=\myroot_2}\MM_2$
admits the presentation\[\langle~A:R~\rangle^{+}\hbox{\quad with\quad}A=A_1\sqcup A_2\hbox{\quad and\quad }R=R_{\selecName_1}\sqcup R_{\selecName_2}\sqcup\{\myruut_1=\myruut_2\},\]
where $\myruut_i\in A_i^*$ is any fixed representative of~$\myroot_i$ for~$i\in\{1,2\}$.

We shall prove that~$\MM$ admits a complemented presentation~$\langle~A:R_\selecName~\rangle^{+}$
where the syntactic right-complement~$\selecName$ essentially extends the syntactic right-complements~$\selecName_1$ and~$\selecName_2$.
Formally, we first simply set
\[\selec{}{\xx}{\yy}=\selec{i}{\xx}{\yy}\hbox{\quad for\quad}(x,y)\in{A}_i^2\hbox{\quad and\quad}i\in\{1,2\}.\]

For~$w\in{A}_i^*$ and~$i\in\{1,2\}$, the root~$\myroot_i$ admits powers which are right-multiples
of the element represented by~$w$,
so one can define~$\mmu_{w}=\min\{m:\selec{i}{\myruut_i^m}{w}=\varepsilon\}$ (see~Figure~\ref{f:amalgam}).
Therefore we set:
\[\selec{}{\xx}{\yy}=\left\{
\begin{array}{ll}
\selec{i}{\xx}{\myruut_i^{\mmu_{\yy}}}&\text{for }\mmu_{\xx}\geq \mmu_{\yy}\text{,}\\
\selec{i}{\xx}{\myruut_i^{\mmu_{\xx}}}\selec{3-i}{\myruut_{3-i}^{\mmu_{\xx}}}{\yy}&\text{for }\mmu_{\xx}<\mmu_{\yy}\text{,}
\end{array}\right.
\hbox{\qquad for~~}
\begin{array}{rcl}
(x,y)	&\in	&{A}_i\times{A}_{3-i}\text{,}\\
i	&\in	&\{1,2\}\text{.}
\end{array}\]
By construction, $R_\selecName$ includes~${R}$.
In particular, for any atom~$\xx_i\in A_i$ left-dividing~$u_i$ for~$i\in\{1,2\}$, we have~$\mmu_{\xx_i}=1$ and~$\selec{}{\xx_i}{\xx_{3-i}}=\selec{i}{\xx_i}{\myruut_{i}}$, hence~$\{\myruut_1=\myruut_2\}\subseteq R_\selecName$.
Conversely, any relation~$x\selec{}{\xx}{\yy}=y\selec{}{\yy}{\xx}$ in~$R_\selecName$
derives from relations in~$R$ .
Indeed, for~$(x,y)\in{A}_i\times{A}_{3-i}$ and say~$\mx<\my$, we have:
\[\begin{array}{rcl}
x\ \selec{}{\xx}{\yy}
&\stack{def}							&x\ \selec{i}{\xx}{\myruut_i^{\mmu_{\xx}}}\ \selec{3-i}{\myruut_{3-i}^{\mmu_{\xx}}}{\yy}\\
&{}_{\phantom{\selecName}}\equiv_{\selecName_i}		&\myruut_i^{\mmu_{\xx}}\ \selec{3-i}{\myruut_{3-i}^{\mmu_{\xx}}}{\yy}\\
&\equiv								&\myruut_{3-i}^{\mmu_{\xx}}\ \selec{3-i}{\myruut_{3-i}^{\mmu_{\xx}}}{\yy}\\
&{}_{\phantom{\selecName_{3-}}}\equiv_{\selecName_{3-i}}				&y\ \selec{3-i}{\yy}{\myruut_{3-i}^{\mmu_{\xx}}}\\
&\stack{def}							&y\ \selec{}{\yy}{\xx}.
\end{array}\]
The case~$\mx=\my$ is even simpler:
\[x\ \selec{}{\xx}{\yy}\
\stack{def}	\ x\ \selec{i}{\xx}{\myruut_i^{\mmu_{\yy}}}
\equiv_{\selecName_i}		\myruut_i^{\mmu_{\xx}}
\equiv					\myruut_{3-i}^{\mmu_{\xx}}
\equiv_{\selecName_{3-i}}		y\ \selec{3-i}{\yy}{\myruut_{3-i}^{\mmu_{\xx}}}
\ \stack{def}\ 					y\ \selec{}{\yy}{\xx}.\]
So $\selecName$ is well-defined and gives to~$\MM$ a right-complemented presentation~$\langle~A:R_\selecName~\rangle^{+}$. 

\medskip
The syntactic right-complement~$\selecName$ is defined on~$A^2$ and, by Lemma~\ref{lem-selec}, it can be uniquely extended by using
\[\selec{}{u}{vw}=\selec{}{u}{v}\ \selec{}{\selec{}{v}{u}}{w}
\hbox{\quad and\quad}\selec{}{vw}{u}=\selec{}{w}{\selec{}{v}{u}}\tag{\ref{selec-extension}}\]
for any~$u,v,w\in A^*$.

\medskip
The point is now to check the $\selecName$-cube condition.
Since $\MM$ is atomic, we only need to check it on~$A$. 
It suffices to take say~$(\xx,\zz)\in A_1^2$ and~$y\in A_2$, since the other cases are either symmetric or trivial.
To make reading easier, $X_i,Y_i$, and~$Z_i$ will denote~$\myruut_i^{\mmu_{\scriptstyle x}}, \myruut_i^{\mmu_{\scriptstyle y}}$,
and~$\myruut_i^{\mmu_{\scriptstyle z}}$ respectively for~$i\in\{1,2\}$.

%\amx
\def\amxx{\hbox{\color{black}{$X_1$}}}
\def\amxy{\hbox{\color{black}{$X_2$}}}
\def\amxz{\hbox{\color{black}{$X_1$}}}
%\amy
\def\amyx{\hbox{\color{black}{$Y_1$}}}
\def\amyy{\hbox{\color{black}{$Y_2$}}}
\def\amyz{\hbox{\color{black}{$Y_1$}}}
%\amz
\def\amzx{\hbox{\color{black}{$Z_1$}}}
\def\amzy{\hbox{\color{black}{$Z_2$}}}
\def\amzz{\hbox{\color{black}{$Z_1$}}}
For $\mx>\mz>\my$, we have
\[
\left\{\begin{array}{rclr}
	\selec{}{\selec{}{\xx}{\yy}}{\selec{}{\xx}{\zz}}
	&=&\selec{1}{\selec{1}{\xx}{\amyx}}{\selec{1}{\xx}{\zz}}\\
	&=&\selec{1}{\selec{1}{\amyx}{\xx}}{\selec{1}{\amyx}{\zz}},										&\footnotesize\text{($\selecName_1$-cube)}\\
	\selec{}{\selec{}{\yy}{\xx}}{\selec{}{\yy}{\zz}}
	&=&\selec{}{\selec{2}{\yy}{\amyy}\ \selec{1}{\amyx}{\xx}}{\selec{2}{\yy}{\amyy}\ \selec{1}{\amyz}{\zz}}		\hspace*{36pt}\\
	&=&\selec{1}{\selec{1}{\amyx}{\xx}}{\selec{1}{\amyx}{\zz}},										&\footnotesize\text{($\selecName$-epsilon)}
\end{array}\right.\]
\[
\left\{\begin{array}{rclr}
	\selec{}{\selec{}{\zz}{\yy}}{\selec{}{\zz}{\xx}}
	&=&\selec{1}{\selec{1}{\zz}{\amyz}}{\selec{1}{\zz}{\xx}}\\
	&=&\selec{1}{\selec{1}{\amyz}{\zz}}{\selec{1}{\amyz}{\xx}},										&\footnotesize\text{($\selecName_1$-cube)}\\
	\selec{}{\selec{}{\yy}{\zz}}{\selec{}{\yy}{\xx}}
	&=&\selec{}{\selec{2}{\yy}{\amyy}\ \selec{1}{\amyz}{\zz}}{\selec{2}{\yy}{\amyy}\ \selec{1}{\amyx}{\xx}}		\hspace*{37pt}\\
	&=&\selec{1}{\selec{1}{\amyz}{\zz}}{\selec{1}{\amyz}{\xx}}.										&\footnotesize\text{($\selecName$-epsilon)}
\end{array}\right.\]
\newpage For $\mx>\my>\mz$, we have
\[
\left\{\begin{array}{rclr}
	\selec{}{\selec{}{\xx}{\yy}}{\selec{}{\xx}{\zz}}
	&=&\selec{1}{\selec{1}{\xx}{\amyx}}{\selec{1}{\xx}{\zz}}\\
	&=&\selec{1}{\selec{1}{\amyx}{\xx}}{\selec{1}{\amyx}{\zz}}										&\footnotesize\text{($\selecName_1$-cube)}\\
	&=&\selec{1}{\selec{1}{\amyx}{\xx}}{\varepsilon}=\varepsilon,									&\footnotesize\text{($\selecName$-epsilon)}\\
	\selec{}{\selec{}{\yy}{\xx}}{\selec{}{\yy}{\zz}}
	&=&\selec{}{\selec{2}{\yy}{\amyy}\ \selec{1}{\amyx}{\xx}}{\selec{2}{\yy}{\amzy}}						\hspace*{67pt}\\
	&=&\selec{}{\selec{1}{\amyx}{\xx}}{\selec{2}{\selec{2}{\yy}{\amyy}}{\selec{2}{\yy}{\amzy}}}				&\footnotesize\text{($\selecName$-extension)}\\
	&=&\selec{}{\selec{1}{\amyx}{\xx}}{\selec{2}{\selec{2}{\amyy}{\yy}}{\selec{2}{\amyy}{\amzy}}}			&\footnotesize\text{($\selecName_2$-cube)}\\
	&=&\selec{}{\selec{1}{\amyx}{\xx}}{\selec{2}{\varepsilon}{\varepsilon}}=\varepsilon,					&\footnotesize\text{($\selecName$-epsilon)}
\end{array}\right.\]
\[
\left\{\begin{array}{rclr}
	\hspace*{0pt}&&\hspace*{-27.8pt}
	\selec{}{\selec{}{\zz}{\yy}}{\selec{}{\zz}{\xx}}
	=\selec{}{\selec{1}{\zz}{\amzz}\ \selec{2}{\amzy}{\yy}}{\selec{1}{\zz}{\xx}}\\
	&=&\selec{}{\selec{2}{\amzy}{\yy}}{\selec{1}{\selec{1}{\zz}{\amzz}}{\selec{1}{\zz}{\xx}}}				&\footnotesize\text{($\selecName$-extension)}\\
	&=&\selec{}{\selec{2}{\amzy}{\yy}}{\selec{1}{\selec{1}{\amzz}{\zz}}{\selec{1}{\amzz}{\xx}}}				&\footnotesize\text{($\selecName_1$-cube)}\\
	&=&\selec{}{\selec{2}{\amzy}{\yy}}{\selec{1}{\amzz}{\xx}}										&\footnotesize\text{($\selecName_1$-epsilon)}\\
	&=&\selec{}{\amzy\ \selec{2}{\amzy}{\yy}}{\xx\ \selec{1}{\xx}{\amzz}}								&\footnotesize\text{($\selecName$-epsilon)}\\
	&=&\selec{}{\amzy\ \selec{2}{\amzy}{\yy}}{\xx}\ \selec{}{\selec{}{\xx}{\amzy\ \selec{2}{\amzy}{\yy}}}{\selec{1}{\xx}{\amzz}}								&\footnotesize\text{($\selecName$-extension)}\\
	&=&\selec{}{\amzy\ \selec{2}{\amzy}{\yy}}{\xx}\ \selec{}{\selec{1}{\xx}{\amzz}\ \selec{}{\selec{1}{\amzz}{\xx}}{\selec{2}{\amzy}{\yy}}}{\selec{1}{\xx}{\amzz}}		&\footnotesize\text{($\selecName$-extension)}\\
	&=&\selec{}{\amzy\ \selec{2}{\amzy}{\yy}}{\xx},												&\footnotesize\text{($\selecName$-epsilon)}\\
	\hspace*{0pt}&&\hspace*{-27.8pt}
	\selec{}{\selec{}{\yy}{\zz}}{\selec{}{\yy}{\xx}}
	=\selec{}{\selec{2}{\yy}{\amzy}}{\selec{2}{\yy}{\amyy}\ \selec{1}{\amyx}{\xx}}\\
	&=&\selec{2}{\selec{2}{\yy}{\amzy}}{\selec{2}{\yy}{\amyy}}\
		\selec{}{\selec{2}{\selec{2}{\yy}{\amyy}}{\selec{2}{\yy}{\amzy}}}{\selec{1}{\amyx}{\xx}} 			&\footnotesize\text{($\selecName$-extension)}\\
	&=&\selec{2}{\selec{2}{\yy}{\amzy}}{\selec{2}{\yy}{\amyy}}\
		\selec{}{\selec{2}{\selec{2}{\amyy}{\yy}}{\selec{2}{\amyy}{\amzy}}}{\selec{1}{\amyx}{\xx}}			&\footnotesize\text{($\selecName_2$-cube)}\\
	&=&\selec{2}{\selec{2}{\yy}{\amzy}}{\selec{2}{\yy}{\amyy}}\ \selec{1}{\amyx}{\xx}						&\footnotesize\text{($\selecName$-epsilon)}\\
	&=&\selec{2}{\selec{2}{\amzy}{\yy}}{\selec{2}{\amzy}{\amyy}}\ \selec{1}{\amyx}{\xx}					&\footnotesize\text{($\selecName_2$-cube)}\\
	&=&\selec{2}{\amzy\ \selec{2}{\amzy}{\yy}}{\amyy}\ \selec{1}{\amyx}{\xx}							&\footnotesize\text{($\selecName_2$-epsilon)}\\
	&=&\selec{}{\amzy\ \selec{2}{\amzy}{\yy}}{\xx}.
\end{array}\right.\]
For $\my>\mx>\mz$, we have
\[
\left\{\begin{array}{rclr}
	\selec{}{\selec{}{\xx}{\yy}}{\selec{}{\xx}{\zz}}												
		&=&\selec{}{\selec{1}{\xx}{\amxx}\ \selec{2}{\amxy}{\yy}}{\selec{1}{\xx}{\zz}}					\hspace*{68pt}\\
		&=&\selec{}{\selec{2}{\amxy}{\yy}}{\selec{1}{\selec{1}{\xx}{\amxx}}{\selec{1}{\xx}{\zz}}}			&\footnotesize\text{($\selecName$-extension)}\\
		&=&\selec{}{\selec{2}{\amxy}{\yy}}{\selec{1}{\selec{1}{\amxx}{\xx}}{\selec{1}{\amxx}{\zz}}}			&\footnotesize\text{($\selecName_1$-cube)}\\
		&=&\selec{2}{\selec{2}{\amxy}{\yy}}{\varepsilon}=\varepsilon,								&\footnotesize\text{($\selecName_1$-epsilon)}\\
	\selec{}{\selec{}{\yy}{\xx}}{\selec{}{\yy}{\zz}}
		&=&\selec{}{\selec{2}{\yy}{\amxy}}{\selec{2}{\yy}{\amzy}}\\
		&=&\selec{}{\selec{2}{\amxy}{\yy}}{\selec{2}{\amxy}{\amzy}}								&\footnotesize\text{($\selecName_2$-cube)}\\
		&=&\selec{2}{\selec{2}{\amxy}{\yy}}{\varepsilon}=\varepsilon,								&\footnotesize\text{($\selecName_2$-epsilon)}
\end{array}\right.\]
\[
\left\{\begin{array}{rclr}
	\selec{}{\selec{}{\zz}{\yy}}{\selec{}{\zz}{\xx}}
	&=&\selec{}{\selec{1}{\zz}{\amzz}\ \selec{2}{\amzy}{\yy}}{\selec{1}{\zz}{\xx}}						\hspace*{72pt}\\
	&=&\selec{}{\selec{2}{\amzy}{\yy}}{\selec{1}{\selec{1}{\zz}{\amzz}}{\selec{1}{\zz}{\xx}}}				&\footnotesize\text{($\selecName$-extension)}\\
	&=&\selec{}{\selec{2}{\amzy}{\yy}}{\selec{1}{\selec{1}{\amzz}{\zz}}{\selec{1}{\amzz}{\xx}}}				&\footnotesize\text{($\selecName_1$-cube)}\\
	&=&\selec{}{\selec{2}{\amzy}{\yy}}{\selec{1}{\varepsilon}{\selec{1}{\amzz}{\xx}}}						&\footnotesize\text{($\selecName_1$-epsilon)}\\	
	&=&\selec{}{\selec{2}{\amzy}{\yy}}{\selec{1}{\amzz}{\xx}}										&\footnotesize\text{($\selecName_1$-epsilon)}\\
	&=&\selec{}{\amzy\ \selec{2}{\amzy}{\yy}}{\xx\ \selec{1}{\xx}{\amzz}}								&\footnotesize\text{($\selecName$-epsilon)}\\
	&=&\selec{2}{\amzy\ \selec{2}{\amzy}{\yy}}{\amxy},\\
	\selec{}{\selec{}{\yy}{\zz}}{\selec{}{\yy}{\xx}}
	&=&\selec{2}{\selec{2}{\yy}{\amzy}}{\selec{2}{\yy}{\amxy}}\\
	&=&\selec{2}{\selec{2}{\amzy}{\yy}}{\selec{2}{\amzy}{\amxy}}									&\footnotesize\text{($\selecName_2$-cube)}\\
	&=&\selec{2}{\amzy\ \selec{2}{\amzy}{\yy}}{\amzy\ \selec{2}{\amzy}{\amxy}}						&\footnotesize\text{($\selecName_2$-epsilon)}\\
	&=&\selec{2}{\amzy\ \selec{2}{\amzy}{\yy}}{\amxy}.\\
\end{array}\right.\]
This completes the proof that the syntactic right-complement associated with~$\langle~A:R_{\selecName}~\rangle^{+}$
satisfies the $\selecName$-cube condition. By Theorem~\ref{thm-criterion},
we deduce that the amalgamated free product~$M=\MM_1\star_{\myroot_1=\myroot_2}\MM_2$
is left-cancellative and admits conditional right-lcms.
A symmetric argument allows to conclude that~$M$ is cancellative and admits both left- and right-conditional lcms.

\medskip
Finally, let~$\D_i$ denote the smallest~Garside element in~$\MM_i$ of which~$\myroot_i$ is a root,
say~$\D_i=\myroot_i^{p_i}$ for~$i\in\{1,2\}$.
We naturally define~$\D=\myroot^{p_1\jR p_2}$ with~$\myroot=\myroot_1=\myroot_2$.
By construction, $\D$ inherits the quasi-centrality from~$\D_1$ and~$\D_2$: for any~$b\in\MM$, there exists an element~$b'\in\MM$ satisfying~$\D b=b'\D$.
Therefore, for any left-divisor~$d$ of~$\D$, say~$\D=db$, 
we find~$\D b=b'\D=b'db$, hence~$\D=b'd$ by right-cancellativity.
Using a symmetric argument, we deduce that the set of its right-divisors coincides with the set of its left-divisors.
The latter includes~$A_1\sqcup A_2$ by definition, thus generates~$\MM$.
Therefore, $\D$ is a Garside element for~$\MM$.
\end{proof}

\begin{rem} An extremal case---when the roots are chosen to correspond to some powers
of minimal Garside elements---has been considered very early in~\cite{dehornoyParis,pPhD}.
\end{rem}

\begin{rem} An amalgamated free product of cancellative monoids~$\MM_1$ and~$\MM_2$
along a (cancellative) monoid~$H$ need not inherit the cancellativity of~$\MM_1$ and~$\MM_2$.
While the unitarity assumption on~$H$ ensures the embedding of both~$\MM_1$ and~$\MM_2$
into the amalgamated free product, it is not sufficient to guarantee the cancellativity of the latter
(see for instance the counterexample of~\cite[Section~3]{howie63}).
\end{rem}

\begin{rem} With different approaches and motivations,
a related result appeared in the context of so-called preGarside monoids~\cite{godelleParis}.
No condition for the existence of a Garside element is considered.
A so-called \emph{special} property on~$H$ (stronger than unitarity) is required.
The latter is far to be satisfied in our cyclic amalgam framework.
\end{rem}

\newcolumntype{E}{!{\vrule width 1.3pt}}

\begin{table}[!b]\scriptsize
\setstretch{1.45}
\begin{center}
\begin{widetable}{\textwidth}{Ec|ccEccc|ccccccc|cE}
\hlinewd{1pt}
\vspace*{-10pt}&&&&&&&&&&&&&\\
	&				&	&		&$\D_{\upkappa}^\frac{1}{3}$&	&				&			&			&$\D_{\upkappa}^\frac{2}{3}$	&				&				&				&\!\!$\D_{\upkappa}^\frac{1}{1}$\!\!\\
\hline
&&\vspace*{-8pt}\hspace*{-22pt}$\hh_{\upkappa}$\hspace*{-25pt}
								&$\tty$	&$\ttx\tty$	&$\tty\ttx$	&$\ttx\tty\ttx\tty$	&$\tty\ttx\tty\ttx$&$\tty\tty$		&$\!\!\ttx\ttx\tty\ttx\tty\!\!$	&$\!\!\ttx\tty\ttx\ttx\tty\!\!$	&$\!\!\tty\ttx\ttx\tty\ttx\!\!$	&$\!\!\tty\ttx\tty\ttx\ttx\!\!$	&$\tty\tty\tty$\\
&\hspace*{-7pt}$\hh_{\upchi}$\hspace*{-12pt}&\hspace*{-35pt}	&		&		&		&				&			&			&					&					&					&					&\\
\hlinewd{1pt}
\vspace*{-10pt}&&&&&&&&&&&&&\\
$\D_{\upchi}^{\frac{1}{4}}$	&$\ttb$\hspace*{-12pt}&		&\!4242\!	&\!4611\!	&\!4611\!	&\!\!\!\!580609\!\!\!\!	&\!\!\!\!580609\!\!\!\!&\!\!\!572541\!\!\!&\!\!\!510474\!\!\!		&\!\!\!510474\!\!\!		&\!\!\!510474\!\!\!		&\!\!\!510474\!\!\!		&\!\!1161\!\!\\		
\vspace*{-10pt}&&&&&&&&&&&&&\\
\hline
\vspace*{-10pt}&&&&&&&&&&&&&\\
$\D_{\upchi}^{\frac{1}{3}}$	&$\tta$\hspace*{-12pt}&		&17		&19		&19		&74				&74			&78			&76					&76					&76					&76					&301\\		
\vspace*{-10pt}&&&&&&&&&&&&&\\
\hline
\vspace*{-10pt}&&&&&&&&&&&&&\\
$\D_{\upchi}^{\frac{1}{2}}$	&$\ttb\ttb$\hspace*{-12pt}&&858		&760		&760		&\!\!12622\!\!		&\!\!12622\!\!	&\!\!11010\!\!	&7904				&7904				&7904				&7904				&89\\
\vspace*{-10pt}&&&&&&&&&&&&&\\
\hline
\vspace*{-10pt}&&&&&&&&&&&&&\\
\vspace*{-7pt}			&$\tta\tta$\hspace*{-12pt}&&71		&70		&70		&249				&249			&253			&206					&206					&206					&206					&831\\
$\D_{\upchi}^{\frac{2}{3}}$	&\vspace*{-7pt}	&&&&&&&&&&&&\\
					&$\ttc\ttb$\hspace*{-12pt}&&57		&57		&57		&161				&161			&157			&140					&140					&140					&140					&501\\	
\vspace*{-10pt}&&&&&&&&&&&&&\\
\hline
\vspace*{-10pt}&&&&&&&&&&&&&\\
$\D_{\upchi}^{\frac{1}{1}}$	\vspace*{-6pt}&$\tta\tta\tta$\hspace*{-12pt}&	&300		&225		&225		&917				&917			&837			&611					&611					&611					&611					&19\\
\vspace*{-13pt}&&&&&&&&&&&&&\\
\vspace*{-1pt}&&&&&&&&&&&&&\\
\hlinewd{1pt}
\end{widetable}
\end{center}
\caption{See~Example~\ref{ex-amalOK}: the number of simples of the Garside monoid~$\MM_{\upchi}\star_{\myroot_\upchi=\myroot_\upkappa}\MM_{\upkappa}$
for roots~$\hh_{\upchi}\in\D_{\upchi}^{\frac{p}{q}}$ with~$\frac{p}{q}\in\left\{\frac{1}{4}, \frac{1}{3}, \frac{1}{2}, \frac{2}{3},
\frac{1}{1}\right\}$
and~$\hh_{\upkappa}\in\D_{\upkappa}^{\frac{p}{q}}$ with~$\frac{p}{q}=\left\{\frac{1}{3},\frac{1}{2},\frac{1}{1}\right\}$.
}
\label{tab-amal}
\end{table}

\begin{exa}\label{e:amalKNOT}The simplest examples are those monoids~$K_{p,q}^{+}=\langle~\ttts,\tttt:\ttts^p=\tttt^q~\rangle^{+}$
obtained with~$\MM_1=\langle~\ttts~\rangle^{+}$, $\MM_2=\langle~\tttt~\rangle^{+}$, and $H=\langle~\tth~\rangle^{+}$
with~$\inject_1:\tth\mapsto\ttts^p$ and~$\inject_2:\tth\mapsto\tttt^q$. They are well-known Garside monoids
associated with torus knot groups (whenever~$p$ and~$q$ are coprime), see~\cite{dehornoyParis,dehornoyTheory,pPhD,pTorus}.
For a generalisation, see Corollary~\ref{c:treeProductGarside} and related Example~\ref{ex-treeKappaOne} below.
\end{exa}

\begin{exa}\label{ex-amalOK} Take again $\Chiral=\langle~\tta,\ttb,\ttc: \tta^2 = \ttb\ttc, \ttb^3 = \ttc\tta~\rangle^{+}$ from~Example~\ref{e:Chiral} and~$\Knuth=\langle~\ttx,\tty:
\ttx\tty\ttx\tty\ttx=\tty\tty~\rangle^{+}$ from~Example~\ref{ex-knuth}. Choose the 
amalgamated submonoid~$H=\langle~\tth~\rangle^{+}$ with~$\inject_1:\tth\mapsto\ttc\ttb\in\D_\upchi^{\frac{2}{3}}$
and~$\inject_2:\tth\mapsto\tty\ttx\tty\ttx\in\D_{\upkappa}^\frac{2}{3}$ for instance.
\begin{figure}[!t]
	\begin{center}
	%\ifdraft	{\tikz{\fill[gray] (0,0) rectangle (8.7,8.7);}}
			{\input{chiknuth.tex}}
	\caption{See Example~\ref{ex-amalOK}: the 161-simple lattice of the amalgamated free product monoid~$\Chiral\star_{\ttc\ttb=\tty\ttx\tty\ttx}\Knuth$.}	\label{fig-amalgam}
	\end{center}
\end{figure}
Then~$\Chiral\star_{\ttc\ttb=\tty\ttx\tty\ttx}\Knuth$ is a Garside monoid with minimal Garside element~$\ttb^{8}=(\ttc\ttb)^3=(\tty\ttx\tty\ttx)^3$ admitting $161$ simples,
whose lattice is displayed in~Figure~\ref{fig-amalgam}.
Table~\ref{tab-amal} further illustrates the wide variety of those Garside monoids obtained as free products
with cyclic amalgamation~$\MM_{\upchi}\star_{\myroot_\upchi=\myroot_\upkappa}\MM_{\upkappa}$,
even when one arbitrarily restricts the roots~$\myroot_{\upchi}$ and~$\myroot_{\upkappa}$
to be chosen respectively from the set~$\D_{\upchi}^{\frac{p}{q}}$ with, say, $\frac{p}{q}\in\left\{\frac{1}{4}, \frac{1}{3}, \frac{1}{2}, \frac{2}{3}, 
\frac{1}{1}\right\}$
and the set~$\D_{\upkappa}^{\frac{p}{q}}$ with~$\frac{p}{q}=\left\{\frac{1}{3},\frac{1}{2},\frac{1}{1}\right\}$.
Such Garside structures may become huge: 
choosing~$\myroot_{\upchi}\in\D_{\upchi}^{\frac{3}{4}}$ for instance, 
we obtain almost 74 million simples for~$\MM_{\upchi}\star_{\tta\ttc=\ttx\ttx\tty\ttx\tty}\MM_{\upkappa}$,
and about seven times for~$\MM_{\upchi}\star_{\ttc\tta=\ttx\ttx\tty\ttx\tty}\MM_{\upkappa}$.
\end{exa}

From Theorem~\ref{thm-amalGarside}, we deduce the following corollary that establishes a complete characterisation of those cyclic amalgamated free products
of Garside groups which are Garside as well.

\begin{corollary}\label{cor-amalGarside}
Let~$\MM_1$ and~$\MM_2$ be some Garside monoids.
The (enveloping group of) the cyclic amalgamated free product~$\MM_1\star_{\myroot_1=\myroot_2}\MM_2$ is Garside
if and only if $\myroot_1$ is a root of some Garside element in~$\MM_1$
and $\myroot_2$ is a root of some Garside element in~$\MM_2$.
\end{corollary}

\begin{proof} Theorem~\ref{thm-amalGarside} coincides with~$(\Leftarrow)$, so it suffices to show~($\Rightarrow$). The centre
 of the amalgamated free product~$\MM=\MM_1\star_{\myroot_1=\myroot_2}\MM_2$ is~$Z(\MM)=Z(\MM_1)\cap Z(\MM_2)$ (see~\cite{MKS} or~\cite{jpp} for instance).
Since the amalgamated submonoid~$H$ is infinite cyclic by hypothesis, $Z(\MM)$ is then either trivial or infinite cyclic.
Now, since~$\MM$ is nontrivial and assumed to be Garside, $Z(\MM)$ cannot be trivial,
hence $Z(\MM)$ is infinite cyclic, say~$\langle~z~\rangle^{+}$.
Therefore, the central Garside element~$z$ is some nontrivial power of~$\hh_i$ for~$i\in\{1,2\}$.
\end{proof}

\medbreak
At this point we have to emphasise the associativity of the free product with cyclic amalgamated submonoids in~Theorem~\ref{thm-amalGarside}. A consequence is the following corollary.

\medbreak
A \emph{weighted tree}~$\TT{~}$ is a tree with vertex set~$\VT{~}$ and edge set~$E(\mathcal{T})$ together with a weight map~$\upomega$ which, with every edge between two vertices~$a$ and~$b$, associates two nonzero integer weights~$\upomega_{a,b}$ and~$\upomega_{b,a}$. Such a double weighting is displayed by using label pairs: $a\edg{11}{\upomega_{a,b}}{\upomega_{b,a}} b$. In this case, the \emph{tree product}~$\GT{~}$ is the group presented by~
\[\langle~a\in\VT{~}:
a^{\upomega_{a,b}}=b^{\upomega_{b,a}}\hbox{~for~}\{a,b\}\in E(\mathcal{T})~\rangle.\]

\begin{corollary}\label{c:treeProductGarside}
Every tree product of infinite cyclic groups is a Garside group.
\end{corollary}

\begin{exa}\label{ex-treeKappaOne} Consider the weighted trees of Figure~\ref{fig-treeKappa}. The \emph{positively reduced} tree~$\TTP{0}$ (right) has vertices with extra gray labels that will be explained later in the proof of~Theorem~\ref{thm-hnnGarside} and in~Example~\ref{ex-treeKappaTwo}. According to Corollary~\ref{c:treeProductGarside}, both trees present the same Garside group~$\GT{0}$.
\begin{figure}[!hb]
\centering
	\begin{tikzpicture}[scale=.75,>=latex,node distance=1.2cm,very thick,top color=white,bottom color=gray]
	\tikzstyle{every state}=[minimum size=12pt,inner sep=1.5pt]
	\node[state] (v1) at (-1.6,-1.1) {\footnotesize\color{gray}$~$};
	\node[state] (v3) at (0.5,-3.9) {\footnotesize\color{gray}$~$};
	\node[state] (v4) at (4.1,.7) {\footnotesize\color{gray}$~$};
	\node[state] (v6) at (2.3,-5) {\footnotesize\color{gray}$~$};
	\node[state] (v7) at (3.5,-3.2) {\footnotesize\color{gray}$~$};
	\node[state] (u1) at (0,0) {\footnotesize\color{gray}$~$};
	\node[state] (u2) at (2,-0.5) {\footnotesize\color{gray}$~$};
	\node[state] (u3) at (0.5,-2) {\footnotesize\color{gray}$~$};
	\node[state] (u4) at (4,-1.5) {\footnotesize\color{gray}$~$};
	\node[state] (u5) at (6,-0.5) {\footnotesize\color{gray}$~$};
	\node[state] (u6) at (5.5,-3.2) {\footnotesize\color{gray}$~$};
	\node[state] (u7) at (4.5,-5) {\footnotesize\color{gray}$~$};
	\path 	(v1) edge node[pos=0.25,above left=-2,scale=.8] {$-1$}
					node[pos=0.9,above left=-1.5,scale=.7] {7} (u1)
	      		(v3) edge node[pos=0.25,left=-2,scale=.8] {$-1$}
					node[pos=0.8,left=-1.5,scale=.7] {191} (u3)
	      		(v4) edge node[pos=0.2,left,scale=.8] {2}
					node[pos=0.8,right=-2,scale=.8] {$-5$} (u4)
	      		(v4) edge node[pos=0.1,above right=-2,scale=.8] {1}
					node[pos=0.8,above right=-2,scale=.8] {2} (u5)
	      		(u1) edge node[pos=0.2,above,scale=.8] {3}
					node[pos=0.8,above=-1,scale=.8] {$-2$} (u2)
	      		(u2) edge node[pos=0.02,below,scale=.8] {4}
					node[pos=0.6,below=1.3,scale=.8] {3} (u3)
	      		(u2) edge node[pos=0.25,above,scale=.8] {3}
					node[pos=0.95,above,scale=.8] {3} (u4)
	      		(u6) edge node[pos=0.2,above right=-2,scale=.8] {2}
					node[pos=0.8,above right=-2,scale=.8] {5} (u4)
			(v7) edge node[pos=0.15,below=-.2,scale=.8] {3}
					node[pos=0.72,below=-.2,scale=.8] {$-3$} (u6)
			(v7) edge node[pos=0.22,above left=-2.2,scale=.8] {$-1$}
					node[pos=0.8,above left=-2,scale=.8] {1} (v6)
			(v6) edge node[pos=0.2,above=-1,scale=.8] {1}
					node[pos=0.75,above=-1,scale=.8] {$-2$} (u7);
	\end{tikzpicture}
	\hspace*{15mm}
	\begin{tikzpicture}[scale=.75,>=latex,node distance=1.2cm,very thick,top color=white,bottom color=gray]
	\tikzstyle{every state}=[minimum size=12pt,inner sep=1.5pt]
	\node[state] (u1) at (0,0) {\footnotesize\color{gray}$90$};
	\node[state] (u2) at (2,-0.5) {\footnotesize\color{gray}$60$};
	\node[state] (u3) at (0.5,-2) {\footnotesize\color{gray}$45$};
	\node[state] (u4) at (4,-1.5) {\footnotesize\color{gray}$60$};
	\node[state] (u5) at (6,-0.5) {\footnotesize\color{gray}$48$};
	\node[state] (u6) at (5.5,-3.2) {\footnotesize\color{gray}$24$};
	\node[state] (u7) at (4.5,-5) {\footnotesize\color{gray}$48$};
	\path 	(u1) edge node[pos=0.2,above,scale=.8] {3}
					node[pos=0.9,above,scale=.8] {2} (u2)
	      		(u2) edge node[pos=0.02,below,scale=.8] {4}
					node[pos=0.6,below=1.3,scale=.8] {3} (u3)
	      		(u2) edge node[pos=0.25,above,scale=.8] {3}
					node[pos=0.95,above,scale=.8] {3} (u4)
	      		(u4) edge node[pos=0.1,above,scale=.8] {5}
					node[pos=0.8,above,scale=.8] {4} (u5)
	      		(u6) edge node[pos=0.01,left,scale=.8] {2}
					node[pos=0.7,left,scale=.8] {5} (u4)
			(u6) edge node[pos=0.3,right=-.2,scale=.8] {3}
					node[pos=0.95,right=-.2,scale=.8] {6} (u7);
	\end{tikzpicture}
\caption{See Examples~\ref{ex-treeKappaOne} and~\ref{ex-treeKappaTwo}: a rather general weighted tree~$\TT{0}$ (left)
and the corresponding \emph{positively reduced} weighted tree~$\TTP{0}$ (right); both generate the same Garside group~$\GT{0}$ by Corollary~\ref{c:treeProductGarside}.
% together with each vertex~$\vv$ labelled with the value of~$\ee(\vv)$
}
\label{fig-treeKappa}
\end{figure}
\end{exa}

\begin{proof}[Proof of Corollary~\ref{c:treeProductGarside}] The point is to show that any tree with nonzero weights can be transformed into another---generating the same group---with weights all belonging to~$\{2,3,4\ldots\}$. Such a transformation actually requires two steps, say a positive transformation and an atomic transformation, which happen to be commutative.

\smallskip The atomic transformation amounts to delete any vertex~$a$ admitting a weight~$\upomega(a,b)=\pm1$ for some vertex~$b$ and, therefore, to connect the latter to each vertex~$c$ among the other possible former neighbours of~$a$ with the weights~$\upomega(c,b)=\upomega(c,a)$ and~$\upomega(b,c)=\pm\upomega(b,a)\upomega(a,c)$. By applying a finite sequence of such so-called elementary collapses, we finally obtain a tree---generating the same group---without weight~$\pm1$, which is well known as a \emph{reduced} weighted tree (see~\cite{levitt} for instance).%While any tree may be reduced by applying a finite sequence of elementary collapses, the reduction is not always unique.

\smallskip Especially relevant for the tree case, the positive transformation ultimately amounts to simply take the absolute value of each weight, but it can be rigorously described and justified as follows. First root the tree at any distinguished vertex~$\vv$. Then apply the following recursive algorithm from~$\vv$ to the leaves:
\begin{itemize}
\item for each vertex~$b\not=\vv$ with parent~$a$, if~$\upomega_{b,a}$ is negative, then take the opposite of each weight~$\upomega_{b,c}$
with~$c$ neighbour of~$b$ (this corresponds with exchanging the generator~$b$ into its inverse~$b^{-1}$);
\item for each vertex~$b$, for each child~$c$ with~$\upomega_{b,c}$ negative, take the opposite of both~$\upomega_{b,c}$ and~$\upomega_{c,b}$
(this corresponds with rewriting the relation~$b^{\upomega_{b,c}}=c^{\upomega_{c,b}}$ into~$b^{-\upomega_{b,c}}=c^{-\upomega_{c,b}}$).
\end{itemize}

Therefore, from any finite weighted tree~$\TT{~}$, we obtain a \emph{positively reduced} weighted tree~$\TTP{~}$---that is, whose weights are all in~$\{2,3,4\ldots\}$---satisfying~$\GT{~}\cong\GTP{~}$. Now, the monoid~$\MTP{~}$ presented by $\langle~a\in\VTP{~}:
a^{\upomega_{a,b}}=b^{\upomega_{b,a}}\hbox{~for~}\{a,b\}\in E(\mathcal{T'})~\rangle^{+}$ is an atomic monoid by Definition~\ref{D:Norm} and thus a Garside monoid by Theorem~\ref{thm-amalGarside} and~Corollary~\ref{cor-amalGarside}, that makes $\GT{~}$ a Garside group.
\end{proof}

%%%%%%%%%%%%%%%%%%%%%%%%%%
%%%%%%%%%%%%%%%%%%%%%%%%%%
%%%%%%%%%%%%%%%%%%%%%%%%%% 
\section{HNN extensions}\label{s:hnnextension}

We characterise and study those cyclic HNN extensions of a Garside monoid whose enveloping groups are Garside as well
(under some mild atomicity condition).
Again the roots of Garside elements play a crucial role.

\medskip
\begin{defin}\label{d:hnnext} Let~$\MM$ and~$H$ be two monoids with morphisms~$\inject_1:H\hookrightarrow\MM$ and~$\inject_2:H\hookrightarrow\MM$.
The \emph{HNN extension} of~$\MM$ with respect to~$H$, $\inject_1$, and~$\inject_2$ is the monoid
\[\langle~\MM,\stable~: \inject_1(\hh)\stable=\stable\inject_2(\hh), \hh\in H~\rangle^{+}.\]
\end{defin}

\begin{theorem}\label{thm-hnnGarside} Let~$\MM$ be a Garside monoid and $H$ be the infinite cyclic monoid~$\langle~\hh~\rangle^{+}$
with a morphism~$\inject_i:H\hookrightarrow\MM$ for~$i\in\{1,2\}$ satisfying~$\length{\inject_1(\hh)}=\length{\inject_2(\hh)}$.
Then the enveloping group of the HNN extension~$\langle~\MM,\stable:\inject_1(\hh)\stable=\stable\inject_2(\hh)~\rangle^{+}$ is a Garside group
if and only if $\inject_1(\hh)$ and~$\inject_2(\hh)$ are two $n$-th roots of the same Garside element in~$\MM$ for some~$n>0$.
\end{theorem}

Here again, a necessary assumption is actually that $\inject_1(H)$ and $\inject_2(H)$ have to contain the same Garside element of~$M$.
When restricted to cyclic HNN extensions, the latter naturally expresses in terms of roots of Garside elements (see also Remark~\ref{rem-extraction}). 

\begin{proof} Let~$\onceMM$ denote the HNN extension~$\langle~\MM,\stable:\hh_1\stable=\stable \hh_2~\rangle^{+}$
where~$\hh_i$ denotes~$\inject_i(\hh)$ for~$i\in\{1,2\}$ and let~$\onceGG$ denote its enveloping group~$\langle~\GG,\stable:\hh_1\stable=\stable \hh_2~\rangle$.

\medskip
By hypothesis, the monoid~$\MM$ is cancellative and the monoid~$H$ is the infinite cyclic monoid,
hence $\inject_i(H)=\langle~\hh_i~\rangle^{+}$ is a so-called \emph{unitary} submonoid of~$\MM$ for each~$i\in\{1,2\}$
(see also the proof of Theorem~\ref{thm-amalGarside}).
Therefore, according to~\cite[Theorem~1]{howie63hnn}, the monoid~$\MM$ embeds into its cyclic HNN extension~$\onceMM$.

\medskip
$(\Rightarrow)$ The centre~$Z(\onceGG)$ of~$\onceGG$ is the subgroup~$Z(\GG)\cap\inject_1(K)\cap\inject_2(K)$ \cite{MKS,jpp}
where~$K$ denotes the group of fractions of~$H$ and the morphism~$\inject_i:K\to\GG$ extends~$\inject_i:H\to\MM$ for each~$i\in\{1,2\}$.
Since~$K$ is infinite cyclic by hypothesis, $Z(\onceGG)$ is then either trivial or infinite cyclic.
Now, since~$\onceGG$ is nontrivial and assumed to be Garside, $Z(\onceGG)$ cannot be trivial,
hence $Z(\onceGG)$ is infinite cyclic.

\medskip
By Corollary~\ref{cor-quasicenter}, $Z(\onceGG) $ is the group of fractions of some~$\langle~z~\rangle^{+}=Z(\MM)\cap\inject_1(H)\cap\inject_2(H)$.
Therefore, the Garside element~$z$ is some nontrivial power~$\hh_1^\ell$ of~$\hh_1$ which has to be central in~$\GG$. 

\medskip
Moreover, $\hh_1^\ell$ has to commute with~$\stable$, that is, $\hh_1^\ell\stable=\stable \hh_1^\ell$ holds.
Now, $\hh_1$ satisfies~$\hh_1\stable=\stable \hh_2$ by definition of~$\onceGG$, hence $\hh_1^d\stable=\stable \hh_2^d$ for any~$d\in\ZZ$.
We find $\stable \hh_1^\ell=\stable \hh_2^\ell$ in~$\onceGG$, hence~$\hh_1^\ell=\hh_2^\ell$ in~$Z(\GG)$.

\medskip 
$(\Leftarrow)$ Let~$A$ be the set of atoms of~$\MM$.
As it stands, the monoid~$\onceMM$ need not provide a Garside structure to its enveloping group~$\onceGG$.
The trick is to introduce the map~$\autom$ defined by
	\[\autom~:~x\mapsto\left\{\begin{array}{ll}
	\stable&\text{~for~}\xx=\stable\text{,}\\
	\xx\stable&\text{~for~}\xx\in A\text{.}\\
	\end{array}\right.\]
%which induces an automorphism of the group~$\onceGG$.
which uniquely extends to a group isomorphism, corresponding to a change of generators.
We claim that the submonoid of~$\autom(\onceMM)$ generated by~$A\sqcup\{\stable\}$ is Garside.
We consider the natural candidate-monoid~$\seccMM$ whose presentation is obtained from that of~$\onceMM$
by keeping the set of atoms~$A\sqcup\{\stable\}$ and by stuffing each relation with the stable atom~$\stable$.
Formally, for any nonempty word~$u=\prod_{i=1}^ma_i$ with~$a_i\in A$,
the word denoted by~$\shuf{\,u\,}{\stable}$ is $a_1\prod_{i=2}^m(\stable a_i)$.
So, every relation~$u=v$ from~$\MM$ becomes~$\shuf{\,u\,}{\stable}=\shuf{\,v\,}{\stable}$ in~$\seccMM$ and the special relation $\hh_1\stable=\stable \hh_2$ becomes
$\shuf{\hh_1}{\stable}\stable=\stable\shuf{\hh_2}{\stable}$ in~$\seccMM$. 
We shall show that~$\seccMM$ is a Garside monoid
whose enveloping group is the HNN extension~$\onceGG$.

\medskip
First, the assumption~$\length{\inject_1(\hh)}=\length{\inject_2(\hh)}$ guarantees~$\seccMM$ to inherit the atomicity from~$\MM$.
Indeed, the mapping~$\seccNu$ from~$\seccMM$ to the integers defined by
\[\seccNu(b_0\,\stable\,b_1\,\stable\,\cdots\,\stable\,b_n)=n+\sum_{k=0}^n\length{b_k}\]
with~$b_k\in\MM$ for~$0\leq k\leq n$ satisfies Condition~(i) from~Lemma~\ref{lem-atomic} (the same holds for~$\onceMM$).

\medskip
Let~$M=\langle~A:R_{\selecName}~\rangle^{+}$. By definition, we have
\[\onceMM=\langle~A\sqcup\{\stable\}:R_{\selecName}\sqcup\{u_1\stable=\stable u_2\}~\rangle^{+},\]
where $\myruut_i\in A^*$ is any fixed representative of~$\hh_i$ for~$i\in\{1,2\}$.
We deduce
	\[\seccMM=\langle~A\sqcup\{\stable\}:\bigsqcup_{(x,y)\in{A}^2}\left\{\shuf{\xx\selec{}{\xx}{\yy}}{\stable}
	=\shuf{\yy\selec{}{\yy}{\xx}}{\stable}\right\}\sqcup\{\shuf{u_1}{\stable}\stable=\stable\shuf{u_2}{\stable}\}~\rangle^{+}.\]
We shall prove that~$\seccMM$ admits a complemented presentation~$\langle~A\sqcup\{\stable\}:R_{\scalebox{.75}{$\seccName$}}~\rangle^{+}$
where the syntactic right-complement~$\seccName$ essentially extends the syntactic right-complement~$\selecName$ via the map~$\autom$.
Formally, we first simply set 
	\[\secc{\xx}{\yy}=\stable\ \shuf{\selec{}{\xx}{\yy}}{\stable}	\hbox{\quad for\quad}	(x,y)\in{A}^2.\]
Let~$x_1\in{A}$ denote the leftmost letter of~$u_1$ and let~$v_1\in{A}^*$ be the word satisfying~$u_1=x_1v_1$. Then we set:
	\[\secc{\xx_1}{\stable}=\stable\ \shuf{v_1}{\stable}\stable
	\hbox{\quad and\quad}
	\secc{\stable}{\xx_1}=\shuf{u_2}{\stable}.\]
Finally, for~$x\in{A}\smallsetminus\{\xx_1\}$, we set:
\[\secc{\xx}{\stable}=
	\stable\ \shuf{\selec{}{\xx}{u_1}}{\stable}
	\hbox{\quad and\quad}
	\secc{\stable}{\xx}=\shuf{u_2}{\stable}\ \shuf{\selec{}{u_1}{\xx}}{\stable}.\]
So $\seccName$ is well-defined and gives to~$\seccMM$
a right-complemented presentation~$\langle~A\sqcup\{\stable\}:R_{\scalebox{.75}{$\seccName$}}~\rangle^{+}$. 

\medskip
The syntactic right-complement~$\seccName$ is defined on~$(A\sqcup\{\stable\})^2$
and, by Lemma~\ref{lem-selec}, it can be uniquely extended by using
\[\secc{u}{vw}=\secc{u}{v}\ \secc{\secc{v}{u}}{w}
\hbox{\quad and\quad}\secc{vw}{u}=\secc{w}{\secc{v}{u}}\tag{$\seccName$-extension}\]
for any~$u,v,w\in (A\sqcup\{\stable\})^*$.
The diagrammatic mechanisms associated with a syntactic complement turn out to be essentially invariant under the stuffing operation.
Indeed, we prove that
\begin{equation}\secc{\shuf{u}{\stable}}{\shuf{v}{\stable}}=\stable\ \shuf{\selec{}{u}{v}}{\stable}	\label{invariance}		\tag{$\seccName$-invariance}		\end{equation}
holds for~$u,v\in{A}^+$.
We use induction on the number~$n$ of steps of computation.
For~$n=1$, hence~$|u|=|v|=1$, we have~$(\shuf{u}{\stable},\shuf{v}{\stable})=(u,v)\in{A}^2$ and the result follows by definition.

Assume~$n>1$ and, without loss of generality, $v=xw$ with~$x\in{A}$ and~$w\not=\varepsilon$,
the $n$ steps decompose into~$n_1+1+n_2$ steps according to the following diagram\nopagebreak[4]
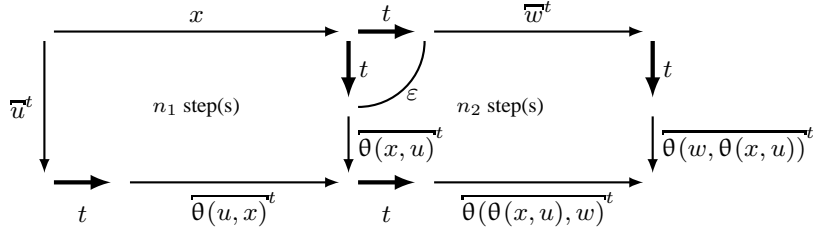
\begin{center}
\begin{tikzpicture}[thick,node distance=10mm]
	\node (00) at (0,2) {};
	\node (01) [right of=00] {};
	\node (02) [right of=01,node distance=30mm] {};
	\node (03) [right of=02] {};
	\node (04) [right of=03,node distance=30mm] {};
	\node (10) [below of=00] {};
	\node (11) [right of=10,node distance=20mm] {\footnotesize $n_1$ step(s)};
	\node (12) [right of=11,node distance=20mm] {};
	\node (13) [right of=12,node distance=20mm] {\footnotesize $n_2$ step(s)};
	\node (14) [right of=13,node distance=20mm] {};
	\node (20) [below of=10] {};
	\node (21) [right of=20] {};
	\node (22) [right of=21,node distance=30mm] {};
	\node (23) [right of=22] {};
	\node (24) [right of=23,node distance=30mm] {};
	\path[->,>=latex]
		(00)	edge 		node[above]{$x$}				(02)
		(03)	edge 		node[above]{$\shuf{w}{\stable}$}	(04)
		(00)	edge 		node[left]{$\shuf{u}{\stable}$}	(20)
		(12)	edge 		node[right]{$\shuf{\selec{}{x}{u}}{\stable}$}			(22)
		(14)	edge 		node[right]{$\shuf{\selec{}{w}{\selec{}{x}{u}}}{\stable}$}			(24)
		(21)	edge 		node[below]{$\shuf{\selec{}{u}{x}}{\stable}$}		(22)
		(23)	edge 		node[below]{$\shuf{\selec{}{\selec{}{x}{u}}{w}}{\stable}$}		(24);
	\path[->,>=latex,ultra thick]
		(02)	edge 		node[above]{$\stable$}			(03)
		(02)	edge 		node[right]{$\stable$}		(12)
		(04)	edge 		node[right]{$\stable$}		(14)
		(20)	edge 		node[below=5]{$\stable$}	(21)
		(22)	edge 		node[below=5]{$\stable$}	(23);
	\path (03) edge[bend left, out=45, in=135] node[below right=-3]{$\varepsilon$} (12);
\end{tikzpicture}
\end{center}
and we obtain
\[\begin{array}{rcl}
\secc{\shuf{v}{\stable}}{\shuf{u}{\stable}}
	&=			&\secc{\shuf{xw}{\stable}}{\shuf{u}{\stable}}\\
	&=			&\secc{x\ \stable\ \shuf{w}{\stable}}{\shuf{u}{\stable}}\\
	&\stack{(IH)}	&\stable\ \shuf{\selec{}{w}{\selec{}{x}{u}}}{\stable}\\
	&=			&\stable\ \shuf{\selec{}{xw}{u}}{\stable}\\
	&=			&\stable\ \shuf{\selec{}{v}{u}}{\stable}\\
\end{array}
\hfill
\begin{array}{rcl}
\quad\secc{\shuf{u}{\stable}}{\shuf{v}{\stable}}
	&=			&\secc{\shuf{u}{\stable}}{\shuf{xw}{\stable}}\\
	&=			&\secc{\shuf{u}{\stable}}{x\ \stable\ \shuf{w}{\stable}}\\
	&\stack{(IH)}	&\stable\ \shuf{\selec{}{u}{x}}{\stable}\ \stable\ \shuf{\selec{}{\selec{}{x}{u}}{w}}{\stable}\\
	&=			&\stable\ \shuf{\selec{}{u}{x}\ \selec{}{\selec{}{x}{u}}{w}}{\stable}\\
	&=			&\stable\ \shuf{\selec{}{u}{v}}{\stable},
\end{array}\]
which concludes the induction. Similarly, we could find
\begin{equation*}
\secc{\shuf{u}{\stable}}{\shuf{v}{\stable}\stable}=\stable\ \shuf{\selec{}{u}{v}}{\stable}
\hbox{\qquad and\qquad}
\secc{\shuf{v}{\stable}\stable}{\shuf{u}{\stable}}=\shuf{\selec{}{v}{u}}{\stable}.	\label{invariance$_2$}		\end{equation*}%\tag{$\seccName$-invariance$_2$}		

\medskip
The point is now to check the $\seccName$-cube condition.
Since $\seccMM$ is atomic, we only need to check it on~$A\sqcup\{\stable\}$. 
It suffices to consider the two following cases, since the others are either symmetric or trivial.

\[
\left\{\begin{array}{rclr}
	\secc{\secc{\xx}{\stable}}{\secc{\xx}{\yy}}
	&=&\secc{\stable\ \shuf{\selec{}{\xx}{u_1}}{\stable}}{\stable\ \shuf{\selec{}{\xx}{\yy}}{\stable}}\\
	&=&\secc{\shuf{\selec{}{\xx}{u_1}}{\stable}}{\shuf{\selec{}{\xx}{\yy}}{\stable}}														&\footnotesize\text{($\seccName$-epsilon)}\\
	&=&\stable\ \shuf{\selec{}{\selec{}{\xx}{u_1}}{\selec{}{\xx}{\yy}}}{\stable}															&\footnotesize\text{($\seccName$-invariance)}\\
	\secc{\secc{\stable}{\xx}}{\secc{\stable}{\yy}}
	&=& \secc{\shuf{u_2}{\stable}\ \shuf{\selec{}{u_1}{\xx}}{\stable}}
			{\shuf{u_2}{\stable}\ \shuf{\selec{}{u_1}{\yy}}{\stable}}							\hspace*{35pt}\\
	&=& \secc{\shuf{\selec{}{u_1}{\xx}}{\stable}}
			{\shuf{\selec{}{u_1}{\yy}}{\stable}}																				&\footnotesize\text{($\seccName$-epsilon)}\\
	&=& \stable\ \shuf{\selec{}	{\selec{}{u_1}{\xx}}
						{\selec{}{u_1}{\yy}}}{\stable}																		&\footnotesize\text{($\selecName$-invariance)}\\
	&=&\stable\ \shuf{\selec{}{\selec{}{\xx}{u_1}}{\selec{}{\xx}{\yy}}}{\stable}															&\footnotesize\text{($\selecName$-cube)}\\
\end{array}\right.\]

\[
\left\{\begin{array}{rclr}
	\secc{\secc{\xx}{\yy}}{\secc{\xx}{\stable}}
	&=&\secc{\stable\ \shuf{\selec{}{\xx}{\yy}}{\stable}}	{\stable\ \shuf{\selec{}{\xx}{u_1}}{\stable}} 	\hspace*{60pt}\\
	&=&\secc{\shuf{\selec{}{\xx}{\yy}}{\stable}}	{\shuf{\selec{}{\xx}{u_1}}{\stable}}														&\footnotesize\text{($\seccName$-epsilon)}\\
	&=&\stable\ \shuf{\selec{}{\selec{}{\xx}{\yy}}{\selec{}{\xx}{u_1}}}{\stable}															&\footnotesize\text{($\seccName$-invariance)}\\
	\secc{\secc{\yy}{\xx}}{\secc{\yy}{\stable}}
	&=&\secc{\stable\ \shuf{\selec{}{\yy}{\xx}}{\stable}}	{\stable\ \shuf{\selec{}{\yy}{u_1}}{\stable}}\\
	&=&\secc{\shuf{\selec{}{\yy}{\xx}}{\stable}}		{\shuf{\selec{}{\yy}{u_1}}{\stable}}													&\footnotesize\text{($\seccName$-epsilon)}\\
	&=&\stable\ \shuf{\selec{}{\selec{}{\yy}{\xx}}{\selec{}{\yy}{u_1}}}{\stable}																&\footnotesize\text{($\seccName$-invariance)}\\
	&=&\stable\ \shuf{\selec{}{\selec{}{\xx}{\yy}}{\selec{}{\xx}{u_1}}}{\stable}																&\footnotesize\text{($\selecName$-cube)}\\
\end{array}\right.\]

According to Theorem~\ref{thm-criterion}, $\seccMM$ is therefore a cancellative monoid admitting conditional lcms.

\medskip Here is the point where the hypothesis on the exponents enters the scene.
The existence of a (central) Garside element is what is needed to conclude that~$\seccMM$
is a Garside monoid whose group of fractions is the enveloping group~$\onceGG$.
Indeed, let~$z$ denote some central Garside element of~$\MM$.
The hypothesis~$\ee(\hh_1)=\ee(\hh_2)$ 
implies that the centre is nontrivial
and then that~$\autom(z)=\shuf{\,z\,}{\stable}\stable$ is a central Garside element in~$\seccMM$,
concluding the proof.
\end{proof}

\begin{rem}\label{rem-patho} The additional assumption~$\length{\inject_1(\hh)}=\length{\inject_2(\hh)}$
can be naturally expressed more generally as the existence of any norm~$\upnu$ satisfying Lemma~\ref{lem-atomic}(i).

Whenever the monoid~$\MM$ admits an additive norm~$\upnu_{+}$,
such an assumption becomes superfluous,
thanks to the required condition on exponents: $\hh_1^{\ell}=\hh_2^{\ell}$
implies~$\upnu_+(\hh_1)=\upnu_+(\hh_2)$.

On the contrary, Example~\ref{ex-hnnKnuth} below %(together with the related Conjecture~\ref{conj-hnnGarside})
illustrates various behaviors
of cyclic HNN extensions of~$\Knuth$
possibly depending on this additional assumption.
\end{rem}

\begin{exa}\label{ex-bs} The simplest examples are those well-known Baumslag--Solitar groups
\[\BS(m,m)=\langle~\ttts,\tttt:\ttts^m\tttt=\tttt\ttts^{m}~\rangle,\]
which are cyclic HNN extension of the infinite cyclic group.
Known as \emph{Baumslag--Solitar monoids} (see~\cite{jacksonBSS,cainBSS,pAutomaticon}),
the associated monoids~$\BS^+(m,m)=\langle~\ttts,\tttt:\ttts^m\tttt=\tttt\ttts^{m}~\rangle^{+}$
happen to be cancellative atomic monoids admitting conditional lcms, but fail to be Garside monoids for~$m>1$.
Now, for any~for~$m>1$, the $\tttt$-stuffed version
(that is, the submonoid of~$\autom(\BS^+(m,m))$ generated by~$\{\ttts,\tttt\}$)
turns out to coincide with a dihedral Artin--Tits monoid
\[\BB^+(\hbox{I}_2(2m))=\langle~\ttts,\tttt:(\ttts\tttt)^{m}=(\tttt\ttts)^m~\rangle^{+}\]
which is Garside and embeds into its group of fractions~$\BS(m,m)$ for~$m>0$.
\end{exa}

\begin{exa}\label{ex-hnnKnuth} Much more complicated examples are provided
by the slightly pathological Garside monoid~$\Knuth=\langle~\ttx,\tty:\ttx\tty\ttx\tty\ttx=\tty\tty~\rangle^{+}$ from Example~\ref{ex-knuth}.
Recall here that $\Knuth$ admits no additive norm, that is, no norm~$\upnu$ 
satisfying~$\upnu(ab)=\upnu(a)+\upnu(b)$ for any~$(a,b)\in\Knuth$.
Its minimal Garside element~$\D=\tty^3$ admits three cube roots 
while~$\D^2$ admits seven:
\[
\D^{\frac{1}{3}}=\{\tty, \ttx\tty, \tty\ttx\}
\quad\hbox{and}\quad
\D^{\frac{2}{3}}=\{\ttx\tty\ttx\tty, \tty\ttx\tty\ttx, \tty\tty, \ttx\ttx\tty\ttx\tty, \tty\ttx\ttx\tty\ttx, \tty\ttx\ttx\tty\ttx, \tty\ttx\tty\ttx\ttx\}.\]
We have~$\length{\tty}=1$, $\length{\ttx\tty}=\length{\tty\ttx}=2$, $\length{\ttx\tty\ttx\tty}=\length{\tty\ttx\tty\ttx}=4$,
$\length{\tty\tty}=\length{\ttx\ttx\tty\ttx\tty}=\length{\tty\ttx\ttx\tty\ttx}=\length{\tty\ttx\ttx\tty\ttx}=\length{\tty\ttx\tty\ttx\ttx}=5$, and~$\length{\tty^3}=6$.

\smallskip
Table~\ref{tab-knuth} gathers the numbers of the possible cyclic HNN extensions of~$\Knuth$
when one arbitrarily restricts a pair~$(\hh_1,\hh_2)$ of images of~$\hh$
to be chosen from the set~$\D^{\frac{p}{q}}$ with~$\frac{p}{q}=\left\{\frac{1}{3},\frac{1}{2},\frac{1}{1}\right\}$.

The nine cases corresponding with the restriction to~$\D^{\frac{1}{3}}$ seem to catch at a glance the whole picture.
For~$(\hh_1,\hh_2)=(g,g)$ with~$g\in\D^{\frac{1}{3}}$, 
	the associated cyclic HNN extensions admit $\ttz$-stuffed presentations
	defining Garside monoids~$\seccMM_{\upkappa,g,g}$ (with respectively~$88,96$, and~$96$ simples).
	The latter two are non-isomorphic, but anti-isomorphic.
The monoid~$\seccMM_{\upkappa,\tty\ttx,\ttx\tty}$ is Garside (with $96$ simples),
is neither isomorphic nor anti-isomorphic to any of the latter two, now it is anti-isomorphic
to the monoid~$\seccMM_{\upkappa,\ttx\tty,\tty\ttx}$ (see Figure~\ref{fig-hnnKnuth96}).
The monoid~$\seccMM_{\upkappa,\ttx\tty,\tty}$ is clearly non-atomic, now its enveloping group 
also envelops the Garside monoid~$\seccMM_{\upkappa,\tty,\ttx\tty}$.
Figure~\ref{fig-hnnKnuth126} displays its $126$-simple lattice, 
whose anti-isomorphic image corresponds to the (Garside) monoid~$\seccMM_{\upkappa,\tty\ttx,\tty}$
(whose enveloping group of fractions also envelops the non-atomic monoid~$\seccMM_{\upkappa,\tty,\tty\ttx}$).

\begin{figure}[!t]
\begin{center}
	\ifdraft	{\tikz{\fill[gray] (0,0) rectangle (5.5,5.5);}}
			{\input{losKnuthHNN1}}
\caption{See~Example~\ref{ex-hnnKnuth}: the (non-isomorphic) 96-simple lattices of~$\seccMM_{\upkappa,\ttx\tty,\ttx\tty}=\langle~\ttx,\tty,\ttz:\ttx\ttz\tty\ttz\ttx\ttz\tty\ttz\ttx=\tty\ttz\tty,\ttx\ttz\tty\ttz=\ttz\ttx\ttz\tty~\rangle^{+}$
(anti-isomorphic to~$\seccMM_{\upkappa,\tty\ttx,\tty\ttx}$) and of~$\seccMM_{\upkappa,\ttx\tty,\tty\ttx}=\langle~\ttx,\tty,\ttz:\ttx\ttz\tty\ttz\ttx\ttz\tty\ttz\ttx=\tty\ttz\tty,\ttx\ttz\tty\ttz=\ttz\tty\ttz\ttx~\rangle^{+}$
(anti-isomorphic to~$\seccMM_{\upkappa,\tty\ttx,\ttx\tty}$).}	\label{fig-hnnKnuth96}
\end{center}
\end{figure}

\begin{figure}[!t]
\begin{center}
\vspace*{-7pt}
	%\ifdraft	{\tikz{\fill[gray] (0,0) rectangle (9,9);}}
			{\input{losKnuthHNN_b_ab}}
\vspace*{-7pt}
\caption{See~Example~\ref{ex-hnnKnuth}: the 126-simple lattices of~$\seccMM_{\upkappa,\tty,\ttx\tty}=\langle~\ttx,\tty,\ttz:\ttx\ttz\tty\ttz\ttx\ttz\tty\ttz\ttx=\tty\ttz\tty,\tty\ttz=\ttz\ttx\ttz\tty~\rangle^{+}$
(anti-isomorphic to~$\seccMM_{\upkappa,\tty\ttx,\tty}$).}	\label{fig-hnnKnuth126}
\end{center}
\end{figure}

\begin{table}[!h]\scriptsize
\setstretch{1.28}
\begin{center}
\begin{widetable}{\textwidth}{Ec|ccEccc|ccccccc|cE}
\hlinewd{1pt}
&&&&&&&&&&&&&\vspace*{-10pt}\\
	&				&	&		&$\D^\frac{1}{3}$&	&				&			&			&$\D^\frac{2}{3}$	&				&				&				&$\D^\frac{1}{1}$\\
\hline
	&		&$\!\!\hh_2\!\!$	&$\tty$	&$\ttx\tty$	&$\tty\ttx$	&$\ttx\tty\ttx\tty$	&$\tty\ttx\tty\ttx$&$\tty\tty$		&$\!\!\ttx\ttx\tty\ttx\tty\!\!$	&$\!\!\ttx\tty\ttx\ttx\tty\!\!$	&$\!\!\tty\ttx\ttx\tty\ttx\!\!$	&$\!\!\tty\ttx\tty\ttx\ttx\!\!$	&$\tty\tty\tty$\\
	&$\hh_1$&$\!\!\!\!\!\length{\hh_i}\!\!\!\!\!$&1	&2		&2		&4				&4			&5			&5				&5				&5				&5				&6\\
\hlinewd{1pt}
\vspace*{-1pt}	&$\tty$	&1	&88		&\ccg126	&\ccg\oo	&				&			&			&				&				&				&				&\\		
\DO	&$\ttx\tty$			&2	&\ccg\oo	&96		&96		&				&			&			&				&				&				&				&\\
	&$\tty\ttx$			&2	&\ccg126	&96		&96		&				&			&			&				&				&				&				&\\
\hline
	&$\ttx\tty\ttx\tty$	&4	&		&		&		&2304			&2304		&\ccg\oo		&\ccg1561		&\ccg1687		&\ccg1561		&\ccg1687		&\\
	&$\tty\ttx\tty\ttx$	&4	&		&		&		&2304			&2304		&\ccg2552	&\ccg1687		&\ccg1561		&\ccg1561		&\ccg\oo			&\\	
\vspace*{-1pt}&$\tty\tty$	&5	&		&		&		&\ccg2552		&\ccg\oo		&1808		&1109			&1109			&1109			&1109			&\\	
\DT	&$\ttx\ttx\tty\ttx\tty$	&5	&		&		&		&\ccg\oo			&\ccg1687	&1109		&804				&804				&804				&804				&\\
	&$\ttx\tty\ttx\ttx\tty$	&5	&		&		&		&\ccg1561		&\ccg1561	&1109		&804				&804				&804				&804				&\\
	&$\tty\ttx\ttx\tty\ttx$	&5	&		&		&		&\ccg1561		&\ccg1687	&1109		&804				&804				&804				&804				&\\
	&$\tty\ttx\tty\ttx\ttx$	&5	&		&		&		&\ccg1687		&\ccg1561	&1109		&804				&804				&804				&804				&\\
\hline
&&&&&&&&&&&&&\vspace*{-11pt}\\
\DD	&$\tty\tty\tty$		&6	&		&		&		&				&			&			&				&				&				&				&44\\
\hlinewd{1pt}
\end{widetable}
\end{center}
\caption{See Example~\ref{ex-hnnKnuth}: the number of simples of the Garside monoid~$\seccMM_{\upkappa}$
for these HNN extensions of~$\Knuth$ with~$\hh_i=\inject_i(\hh)\in\ \,\DO\,\ \sqcup\ \,\DT\,\ \sqcup\ \,\DD\,\ $ for~$i\in\{1,2\}$. 
Gray cells correspond to pairs~$(\hh_1,\hh_2)$ with~$\length{\hh_1}\not=\length{\hh_2}$,
while the symbol~\oo\ indicates when $\seccMM_{\upkappa}$ is effectively non-atomic.}
\label{tab-knuth}
\end{table}

\end{exa}

The previous result can be fully applied to the class of tree products of infinite cyclic groups.
Gathering~Corollary~\ref{c:treeProductGarside} and Theorem~\ref{thm-hnnGarside}, we obtain:

\begin{corollary}\label{c:hnnTreeProductGarside} Let~$\TT{~}$ be a finite weighted tree. For any vertices~$\vvx,\vvz\in\VT{~}$,
the cyclic HNN extension~$\langle~\GT{~},\stable:\stable\vvx=\vvz\stable~\rangle$ is Garside %a Garside group
%if and only if the weight~$(p_1,q_1,\ldots,p_m,q_m)$ of the path between~$\vvx$ and~$\vvz$ satisfies~$p_1\cdots p_m=q_1\cdots q_m$.
if and only if the weighted path~$\vvx\edg{6}{p_1}{q_1}\cdots\edg{6}{p_m}{q_m}z$ satisfies~$p_1\cdots p_m=q_1\cdots q_m$.
\end{corollary}

%EXPONENTS BY SUCCESSIVE TRIANGULATIONS
Let us first state a straightforward fact about exponents.

\begin{lemma}\label{l:KeyExponent} Let~$\TT{~}$ be a finite weighted tree with weights in~$\{2,3,4\ldots\}$.
For any adjacent vertices~$\vvx,\vvz\in\VT{~}$ with~$\vvx\edg{4.5}{p}{q}\vvz$,
the exponents of~$\vvx$ and~$\vvz$ satisfy~$\frac{\ee(\vvx)}{p}=\frac{\ee(\vvz)}{q}$.
\end{lemma}

\begin{proof}Recall first that, under the mild hypothesis that the weights all belong to~$\{2,3,4\ldots\}$, the monoid~$\MT{~}$ is Garside by (the proof of) Corollary~\ref{c:treeProductGarside}. Then, by definition of~$\ee$ and by hypothesis, we have\[(\vvx^{\ee(\vvx)})^q\stack{def}(\vvz^{\ee(\vvz)})^q=\left(\vvz^q\right)^{\ee(\vvz)}\stack{hyp}\left(\vvx^p\right)^{\ee(\vvz)}.\]We conclude $\frac{\ee(\vvx)}{p}=\frac{\ee(\vvz)}{q}$ by invoking the torsion-freeness of any Garside monoid.
\end{proof}

There are many equivalent ways to compute the exponents of the vertices. We can start from the initial positively reduced tree~$\TT{~}$ and complete it into the associated weighted complete graph~$\hbox{K}_{\TT{~}}$ by successive triangulations: whenever $x\edg{4.5}{a}{b}y$ and $y\edg{4.5}{c}{d}z$ are edges, we add the edge
\[x\edg{12}{a\frac{c\jR b}{b}}{d\frac{b\jR c}{c}}z\]
if needed, according to Lemma~\ref{l:KeyExponent}. The exponent~$\ee(\vv)$ of any vertex~$\vv$ is then obtained by taking the lcm of all these weights~$p$ with~$\vv\edg{5}{p}{q}\vv'$ an edge in~$\hbox{K}_{\TT{~}}$.

% {\color{red} Note that the monoid~$\MT{~}$ is a quotient of the monoid~$\MTP{~}$ for any maximal subtree~$\TTP{~}$ of~$\hbox{K}_{\TT{~}}$.}
%\[\ee(v)=\bigvee_{v \scalebox{.7}{\edg{3}{p}{\cdot}} v'\in\,\hbox{K}_{\TT{~}}}p.\]

\medskip Instead of completing~$\TT{~}$ into~$\hbox{K}_{\TT{~}}$, we can also choose to compute the exponent of any vertex~$\vv$ by rooting the tree~${\TT{~}}_{\!\!\vv}$ at~$\vv$ and by exploring the obtained rooted tree~${\TT{~}}_{\!\!\vv}$.
We have:

\begin{lemma}\label{l:TreeProductExponent} Let~$\TT{~}$ be a finite weighted tree with weights in~$\{2,3,4\ldots\}$.
The exponent~$\ee(\vv)$ %in~$\MT{~}$
of a vertex~$\vv\in\VT{~}$ can be obtained as~$\uplambda({\TT{~}}_{\!\!\vv})$ where~$\uplambda$ is recursively defined by~$\uplambda\left(\leaf{\vvx}\right)=1$ for~$\vvx$ a leaf
and by\[\uplambda\left(\rootree{\vv}{9}{p_1}{q_1}{p_m}{q_m}{\vv_1}{\vv_m}\right)=\bigvee_{1\leq i\leq m}p_i\frac{\uplambda\left(\rootreeS{\vv_i}\right)\jR q_i}{q_i}\] otherwise.
\end{lemma}

\begin{proof} The %very
point is that, for any vertex~$\vvx\in\VT{~}$, the value of~$\uplambda\left(\rootreeS{\vvx}\right)$ coincides with the exponent of the vertex~$\vvx$ in the %(parabolic)
submonoid generated by the subtree pending from this root~$\vvx$.
Up to duplicating some vertices, we can suppose that~${\TT{~}}_{\!\!\vv}$ is either a linear or a starlike tree with root~$\vv$. We use an induction on the number~$m\geq 1$ of branches of this starlike tree.
For~$m=1$, we use an induction on the length~$h\geq 1$
%of the linear rooted tree~$\vv=\vv_0\edgo{6}{p_1}{q_1}\vv_1\cdots\vv_{h-1}\edgo{6}{p_h}{q_h}\vv_h$.
of the linear rooted tree~$\rootreeH{\vv}{9}{p_1}{q_1}{\vv_1}\cdots\rootreeH{\phantom{\vv_{h}}}{9}{p_h}{q_h}{\vv_h}$.
For~$h=1$, we find~$\displaystyle\ee(\vv)=\uplambda\left(\rootreeH{\vv}{9}{p_1}{q_1}{\vv_1}\right)=p_1\frac{\uplambda\left(\leaf{\vv_1}\right)\jR q_1}{q_1}=p_1\frac{1\jR q_1}{q_1}=p_1$ as expected. 
Assume~$h>1$, the exponent of the vertex~$\vv_1$ in the %(parabolic)
submonoid generated by the subtree pending from~$\vv_1$ is $\uplambda\left(\rootreeS{\vv_1}\right)$ by induction hypothesis. We obtain
\[\displaystyle\ee(\vv)=\uplambda\left(\rootreeS{\vv}\right)=p_1\frac{\uplambda\left(\rootreeS{\vv_1}\right)\jR q_1}{q_1},\]
which concludes the induction on~$h\geq 1$.

\medskip Assume~$m>1$ and let~$f_i=p_i\frac{\uplambda\left(\rootreeS{\vv_i}\right)\jR q_i}{q_i}$ for~$1\leq i\leq m$. The exponent of~$\vv$ in~$\MT{~}$ then coincides with the exponent of~$\vv'$ in the monoid~$\MTP{~}$ generated by the star~$\TTP{~}$ below:
\[\staroot{\vv'_{\phantom{\!n\!}}}{11}{f_1}{f_1'}{f_m}{f_m'}{\vv'_1}{\vv'_m}\]
which concludes the induction on~$m\geq 1$ and the proof.
\end{proof}

\begin{proof}[Proof of Corollary~\ref{c:hnnTreeProductGarside}]
By Corollary~\ref{c:treeProductGarside}, $\GT{~}$ is indeed a Garside group.
According to Theorem~\ref{thm-hnnGarside}, it remains to verify that the two vertices~$\vvx,\vvz$ satisfy~$\ee(\vvx)=\ee(\vvz)$ if and only if the weighted path~$\vvx=\vv_0\edg{8}{p_1}{q_1}\vv_1\cdots\vv_{m-1}\edg{8}{p_m}{q_m}\vv_m=\vvz$ satisfies~$p_1\cdots p_m=q_1\cdots q_m$. Note that the latter equality is invariant under both atomic and positive transformations,
that we defined for the proof of Corollary~\ref{c:treeProductGarside}. So we can assume from now on that the weights of the considered tree~$\TT{~}$ all belong to~$\{2,3,\ldots\}$.

\medskip Rooting~$\TT{~}$ at~$\vvx=\vv_0$, we consider the residual subtrees~$\TR{\vv_i}$ rooted at each vertex~$\vv_i$ for~$0\leq i\leq m$ according to the following combing
\begin{center}
	\begin{tikzpicture}[top color=white,bottom color=gray]
		\path[semithick] (0,0) edge[->-=.84,>=latex]	node[pos=0.3,above,scale=.7] {$p_1$}			%->-=.84,>=latex
											node[pos=0.7,above,scale=.7] {$q_1$} (11ex,0);
		\path[semithick] (11ex,0) edge[->-=.84,>=latex]	node[pos=0.3,above,scale=.7] {$p_2$}			%->-=.84,>=latex
											node[pos=0.7,above,scale=.7] {$q_2$} (22ex,0);
		\path[semithick] (22ex,0) edge[->-=.84,>=latex,dotted]						    (33ex,0);	%->-=.84,>=latex,
		\path[semithick] (33ex,0) edge[->-=.84,>=latex]	node[pos=0.3,above,scale=.7] {$p_m$}			%->-=.84,>=latex
											node[pos=0.7,above,scale=.7] {$q_m$} (44ex,0);
		\fill[left color=lightgray,right color=white] (00ex+.25ex,-.25ex) -- (00ex+7ex,-3ex) -- (00ex+5ex,-6ex) -- cycle;
		\node [scale=.85] at (00ex+1.5ex,-4ex)	{$\TR{\vv_0}$};
		\fill[left color=lightgray,right color=white] (11ex+.25ex,-.25ex) -- (11ex+7ex,-3ex) -- (11ex+5ex,-6ex) -- cycle;
		\node [scale=.85] at (11ex+1.5ex,-4ex)	{$\TR{\vv_1}$};
		\fill[left color=lightgray,right color=white] (22ex+.25ex,-.25ex) -- (22ex+7ex,-3ex) -- (22ex+5ex,-6ex) -- cycle;
		\node [scale=.85] at (22ex+1.5ex,-4ex)	{$\TR{\vv_2}$};
		\fill[left color=lightgray,right color=white] (33ex+.25ex,-.25ex) -- (33ex+7ex,-3ex) -- (33ex+5ex,-6ex) -- cycle;
		\fill[left color=lightgray,right color=white] (44ex+.25ex,-.25ex) -- (44ex+7ex,-3ex) -- (44ex+5ex,-6ex) -- cycle;
		\node [scale=.85] at (44ex+1.5ex,-4ex)	{$\TR{\vv_m}$};
		\node [circle,draw=black,fill=lightgray,inner sep=1pt,minimum size=19pt,scale=.85] at (0,0)	{$\vv_0$};
		\node [circle,draw=black,fill=lightgray,inner sep=1pt,minimum size=19pt,scale=.85] at (11ex,0)	{$\vv_1$};
		\node [circle,draw=black,fill=lightgray,inner sep=1pt,minimum size=19pt,scale=.85] at (22ex,0)	{$\vv_2$};
		\node [circle,draw=black,fill=lightgray,inner sep=1pt,minimum size=19pt,scale=.85] at (33ex,0)	{\phantom{$\vv_m$}};
		\node [circle,draw=black,fill=lightgray,inner sep=1pt,minimum size=19pt,scale=.85] at (44ex,0)	{$\vv_m$};
	\end{tikzpicture}
\end{center}
that is, the subtrees~$\TR{\vv_i}$ that we would obtain after erasing from~$\TT{~}$ each edge between~$v_{i-1}$ and~$v_{i}$ for~$1\leq i\leq m$. Note that the definition of these rooted subtrees~$\TR{\vv_i}$ does not depend on the choice of the original root for~$\TT{~}$, provided that such a root is chosen among the vertices~$\vv_i$ along the path between~$\vvx$ and~$\vvz$. We shall compare~$\ee(\vv_0)$ and~$\ee(\vv_m)$ by using the auxiliary map~$\uplambda$ and~Lemma~\ref{l:TreeProductExponent}.

\medskip Let~$r_i$ denote $\uplambda\left(\rootreeR{\vv_i}\right)$ for~$0\leq i\leq m$ and~$e_i$ denote~$\uplambda\left(\rootreeRS{\vv_i}\right)$ for~$0\leq i<m$.
We now show $e_0=r_0\vee \frac{1}{q_1\cdots q_m}\bigvee_{1\leq i\leq m}p_1\cdots p_{i-1}p_i(q_i\vee r_i)q_{i+1}\cdots q_m$ by induction on~$m\geq 0$.
The result is obvious for~$m=0$. Assume~$m>0$. We find first $\displaystyle\ee(\vv_0)=e_{0}=r_0\vee p_1\frac{q_1\vee e_1}{q_1}$ by~Lemma~\ref{l:TreeProductExponent} and then
\begin{align*}
%e_{0}	&\stackrel{\text{\tiny (Lem.~\ref{l:TreeProductExponent})}}{=} 	r_0\vee p_1\frac{q_1\vee e_1}{q_1}\\
e_{0}	&\stackrel{\text{\tiny(IH)}}{=} 			r_0\vee p_1\frac{q_1\vee \left(r_1\vee \frac{1}{q_2\cdots q_m}\bigvee_{2\leq i\leq m}p_2\cdots p_{i-1}p_i(q_i\vee r_i)q_{i+1}\cdots q_m\right)}{q_1}\\
		&\stackrel{\text{\phantom{\tiny(IH)}}}{=} 	r_0\vee\frac{p_1}{q_1q_2\cdots q_{m}}\left((q_1\vee r_1)q_2\cdots q_{m}\vee \bigvee_{2\leq i\leq m}p_2\cdots p_{i-1}p_i(q_i\vee r_i)q_{i+1}\cdots q_m\right)\\
		&\stackrel{\text{\phantom{\tiny(IH)}}}{=} 	r_0\vee\frac{1}{q_1q_2\cdots q_{m}}\left(p_1(q_1\vee r_1)q_2\cdots q_{m}\vee \bigvee_{2\leq i\leq m}p_1p_2\cdots p_{i-1}p_i(q_i\vee r_i)q_{i+1}\cdots q_m\right)\\
		&\stackrel{\text{\phantom{\tiny(IH)}}}{=} 	r_0\vee\frac{1}{q_1q_2\cdots q_{m}}\bigvee_{1\leq i\leq m}p_1\cdots p_{i-1}p_i(q_i\vee r_i)q_{i+1}\cdots q_m,
\end{align*}
which concludes the induction. By just introducing a dummy weight~$q_0=1$, we obtain the completely homogeneous formula
\[\ee(\vvx)=\ee(\vv_0)=\frac{1}{q_1q_2\cdots q_{m}}\bigvee_{0\leq i\leq m}p_1\cdots p_i(q_i\vee r_i)q_{i+1}\cdots q_m.\]
We obtain symmetrically
\[\ee(\vvz)=\ee(\vv_m)=\frac{1}{p_1p_2\cdots p _{m}}\bigvee_{0\leq i\leq m}q_m\cdots q_{i+1}(p_{i+1}\vee r_i)p_i\cdots p_1,\]
where~$p_{m+1}=1$ is also used as some dummy weight.
By using the symmetry of the just obtained iterative versions, we deduce
\[\ee(\vv_0)=\ee(\vv_m)\iff p_1\cdots p_{m}=q_1\cdots q_{m},\]
which concludes the proof.
\end{proof}

\begin{exa}\label{ex-treeKappaTwo} In Figure~\ref{fig-treeKappa}, each vertex~$\vv\in\VTP{0}$ has been labelled with the value~$\ee(\vv)$. So, in addition to the seven direct HNN extensions corresponding to loops (that is, with~$\vvx=\vvz$), Theorem~\ref{thm-hnnGarside} and Corollary~\ref{c:hnnTreeProductGarside} allow to foresee that exactly two single cyclic HNN extensions are Garside (namely, $\ee(\vvx)=\ee(\vvz)\in\{48,60\}$). Multiple cyclic HNN extensions (for instance, $\ee(\vvx)=\ee(\vvz)=48$ plus~$\ee(\vvx')=\ee(\vvz')=60$) provide again Garside groups.
\end{exa}

All the arguments advanced here can be carefully revisited in order to obtain
an even more general result that a Generalised Baumslag--Solitar group is Garside
if and only if its centre is nontrivial,
which essentially seems that a GBS group is Garside whenever it is a GBS-tree group~\cite{levitt}.
We refer to~\cite{DelgadoRobinsonTimm} for a recent, independent and global approach of the calculus of the centre of GBS groups (see also~\cite{Robinson}).

\medskip
Finally, the possible mechanism---if~$\seccMM_{\upkappa,\hh_1,\hh_2}$ is not atomic,
then $\seccMM_{\upkappa,\hh_2,\hh_1}$ is atomic---observed for instance in Example~\ref{ex-hnnKnuth}
and~Table~\ref{tab-knuth},
could lead to raise the question of weakening or even skipping the length assumption in Theorem~\ref{thm-hnnGarside}.
However, as already mentioned in Example~\ref{ex-bs},
the Baumslag--Solitar group~$\BS(m,n)=\langle~\ttts,\tttt:\ttts^m\tttt=\tttt\ttts^{n}~\rangle$
is a Garside group if and only if $\BS(m,n)$ is an automatic group if and only if $m=n$ holds
(see~\cite[Example~7.4.1]{EpsteinWord}).

%\begin{conjecture}\label{conj-hnnGarside} Let~$\MM$ be a Garside monoid and $H$ be the infinite cyclic monoid~$\langle~\hh~\rangle^{+}$
%with a morphism~$\inject_i:H\hookrightarrow\MM$ for~$i\in\{1,2\}$.
%Then the enveloping group of the HNN extension~$\langle~\MM,\stable:\inject_1(\hh)\stable=\stable\inject_2(\hh)~\rangle^{+}$ is a Garside group
%if and only if $\inject_1(\hh)$ and~$\inject_2(\hh)$ are two $n$-th roots of a same Garside element in~$\MM$ for a same~$n>0$.
%\end{conjecture}

%%%%%%%%%%%%%%%%%%%%%%%%%%
%%%%%%%%%%%%%%%%%%%%%%%%%%
%%%%%%%%%%%%%%%%%%%%%%%%%% 
\section{Pietrowski groups}\label{s:pietrowski}

Using the solution of~Pietrowski for the isomorphism problem for one-relator groups with nontrivial centre~\cite{pietrowski},
we finally state that a non-cyclic one-relator group is Garside if and only if its centre is nontrivial.

\begin{theorem}\label{p:pietrowski}{\rm \cite[Theorems~1 and~3]{pietrowski}} Assume that~$G$ is a non-cyclic
one-relator group with nontrivial centre. Then, if $G/[G,G]$ is not free abelian,
$G$ can be uniquely\footnote{modulo mirror symmetry} presented as
\[\langle~\tta_1,\tta_2,\ldots,\tta_m:
\tta_1^{p_1}=\tta_2^{q_1},\ldots,\tta_{m-1}^{p_{m-1}}=\tta_m^{q_{m-1}}~\rangle\]
with~$p_i,q_i\geq 2$ and~$p_i\wedge q_j=1$ for~$i>j$, and, if $G/[G,G]$ is free
abelian, $G$ can be uniquely\footnote{modulo exchange
and cyclic permutations} presented as
\[\langle~\tta,\tta_1,\tta_2,\ldots,\tta_m:
\tta\tta_1=\tta_m\tta,
\tta_1^{p_1}=\tta_2^{q_1},\ldots,\tta_{m-1}^{p_{m-1}}=\tta_{m\phantom{-1}}^{q_{m-1}}~\rangle\]
with~$p_i,q_i\geq 2$, $p_i\wedge q_j=1$ for~$i>j$, and~$p_1p_2\cdots
p_{m-1}=q_1q_2\cdots q_{m-1}$.
\end{theorem}

\begin{defin}\label{d:pietrowski}A group is called \emph{an
$m$-Pietrowski group} if it admits (exactly) one of the presentations of
Theorem~\ref{p:pietrowski} for some (unique) integer~$m>1$.
\end{defin}

The just-above term can be viewed as a short for what we could call a
Baumslag--Collins--Karrass--McCool--Meskin--Metaftsis--Magnus--Murasugi--Pietrowski--Solitar--Steinberg--Taylor group.

\medbreak We gather three properties (see~\cite{meskinPietrowskiSteinberg,Collins1977,McCool1991})
illustrating the still mysterious distribution of one-relator groups with nontrivial centre within the class of Pietrowski groups.

\begin{itemize}
\item[$\bullet$] Every $m$-Pietrowski group with~$m\leq 3$ is a one-relator group: the group
\[\langle~\tta_1,\tta_2,\tta_3:\tta_1^{p_1}=\tta_2^{q_1},\tta_2^{p_{2}}=\tta_3^{q_2}~\rangle\]
with~\(p_i,q_i\geq2\) and~\(q_1\wedge p_2=1\) is isomorphic to the group
\[\langle~\tta_1,\tta_3:\pi_{p_2,q_1}(\tta_1^{p_1},\tta_{3}^{-q_2})~\rangle,\]
where~$\pi_{\lambda,\mu}(x,y)$ denotes the unique primitive element (up to conjugation) with exponent sum~$\lambda$ on~$x$ and~$\mu$ on~$y$
(note that there is a misprint in the original statement of~\cite{meskinPietrowskiSteinberg}).

\item[$\bullet$] Any group~$\langle~\tta_1,\tta_2,\tta_3,\tta_4:\tta_1^2=\tta_2^2,\tta_2^{p_2}=\tta_3^{p_2},\tta_3^3=\tta_4^3~\rangle$ with~$p_2\geq 2$
happens to not be a one-relator group.

\item[$\bullet$] For every~$m$, there exists an $m$-Pietrowski group which is a one-relator group.
\end{itemize}

\bigbreak Applying Corollaries~\ref{c:treeProductGarside} and~\ref{c:hnnTreeProductGarside}, we deduce these rather amazing facts:

\begin{corollary}Every Pietrowski group is Garside.
\end{corollary}

\begin{corollary}\label{cor-iff}A non-cyclic one-relator group is Garside if and only if its centre is nontrivial.
\end{corollary}

Note that Baumslag and Taylor \cite{bt} have given an algorithm for deciding whether or not a one-relator group has a nontrivial centre.
By Corollary~\ref{cor-iff}, one can decide whether or not a given one-relator group is a Garside group.

\begin{exa}\label{ex-oneRelator} Figure~\ref{fig-pietrow} displays the 1353-simple lattice of the Garside structure for the one-relator group
\[\langle~\tta,\ttx:\ttx^8\tta\ttx^{-6}\tta^{-1}\ttx^4\tta\ttx^{-6}\tta^{-1}~\rangle,\]
whose associated Pietrowski presentation is
\[\langle~\tta,\ttx_1,\ttx_2,\ttx_3:\tta\ttx_1=\ttx_3\tta,\ttx_1^4=\ttx_2^2,\ttx_2^3=\ttx_3^6~\rangle.\]
\end{exa}

\begin{figure}[!ht]
	\begin{center}
	\ifdraft	{\tikz{\fill[gray] (0,0) rectangle (10.25,10.25);}}
			{\input{losPietrow4236x13.tex}}
	\caption{The 1353-simple lattice of the Garside structure for the one-relator group~$\langle~\tta,\ttx:\ttx^8\tta\ttx^{-6}\tta^{-1}\ttx^4\tta\ttx^{-6}\tta^{-1}\rangle$ from~Example~\ref{ex-oneRelator}.}\label{fig-pietrow}
	\end{center}
\end{figure}

Interestingly, such a Garside structure provides an explicit biautomatic structure~\cite{dehornoyGarside,pTorus,pTransducer,myhdr}.

\begin{corollary}\label{cor-conj} All Pietrowski groups and, in particular, all (non-cyclic) one-relator groups with nontrivial centre are biautomatic and have solvable conjugacy problem.
\end{corollary}

All Pietrowski groups and, in particular, all (non-cyclic) one-relator groups with nontrivial centre are torsion-free.
In particular, gathering Corollary~\ref{cor-conj} and the solution of Newman for one-relator groups with torsion \cite{newman},
the conjugacy problem remains open for torsion-free one-relator groups with trivial centre (see~\cite[Problem~(O5)]{BMSopen} and also~\cite{LS}).

%%%%%%%%%%%%%%%%%%%%%%%%%%
%%%%%%%%%%%%%%%%%%%%%%%%%%
%%%%%%%%%%%%%%%%%%%%%%%%%% 
%%\bigskip
%%\footnotesize
%%\noindent\textit{Acknowledgments.}
%%
%%The author wishes to acknowledge and express his gratitude to Gilbert~Levitt for his curiosity and to Patrick~Dehornoy
%%for his challenging and friendly support.
%%%The author also thanks the anonymous referee, whose comments greatly improved the paper.

\nocite{serreAmalgames,higmanNeumannNeumann,myhdr}
\bibliographystyle{plain}
\bibliography{garside}

%
%
%
%\bigskip
%\footnotesize
%\noindent\textit{Acknowledgments.}
%This research was partly supported by NSF (grant no. XXXX).
%
%\begin{thebibliography}{SK}
%
%%% Use the widest label as parameter.
%
%%% Reference items may be numbered or have labels of your choice.
%%% The author's surname PRECEDES the initial of the first name
%%% The issue number is only given when the issues are paginated separately.
%%% In book titles, first letters are capitalized.
%%% Only journal volume numbers are boldfaced.
%
%%%%%%%%%%%% To ease editing, use normal size:
%
%\normalsize
%\baselineskip=17pt
%
%%%%%%%%%%%%%%
%
%\bibitem[B]{Barlow} 
%Barlow, M. T.: 
%Diffusions on fractals.  
%In: Lectures on Probability Theory and Statistics (Saint-Flour, 1995), 
%Lecture Notes in Math. 1690, Springer, Berlin, 1--121 (1998)
%
%
%\bibitem[G]{Gratzer} 
%Gr\"atzer, G.:
%More Math into \LaTeX.
%4th ed., Springer, Berlin (2007)
%
%
%\bibitem[SK]{SatoKashiwaraKawai}
%Sato, M., Kashiwara, M., Kawai, M.: 
%Linear differential equations of infinite order and theta functions.
%Adv. Math. \textbf{47}, 300--325 (1983)
%
%
%\bibitem[Sh]{Shchepin}
%Shchepin, E. V.:
%On mappings of the two-dimensional sphere.  
%Uspekhi Mat. Nauk \textbf{58}, no.~2, 169--170 (2003) (in Russian); 
%English transl.: Russian Math. Surveys \textbf{58}, 1218--1219 (2003)
%
%\bibitem[S]{Smith} 
%Smith, J.:
%A new upper bound on the cross number.
%arXiv:2056.7895.
%
%
%\bibitem[V]{Verkaar}
%Verkaar, M.:
%Continuous local martingales and stochastic integration in Banach spaces.
%PhD thesis, Univ. of Helsinki (2001)
%\end{thebibliography}

\end{document}